\documentclass[12pt,a4paper]{article}%
\usepackage{amsmath}
\usepackage{amsfonts}
\usepackage{amssymb}
\usepackage[pdftex]{graphicx}
\usepackage{geometry}
\usepackage{setspace}%
\providecommand{\U}[1]{\protect\rule{.1in}{.1in}}
\newtheorem{theorem}{Theorem}

\newtheorem{remark}[theorem]{Remark}

\begin{document}

\title{Comparative Study of Homotopy Analysis and Renormalization Group Methods on
Rayleigh and\\Van der Pol Equations }
\author{Aniruddha Palit$^{1}$\thanks{ email:mail2apalit@gmail.com} and Dhurjati Prasad
Datta$^{2}$\thanks{Corresponding author; email:dp${_{-}}$datta@yahoo.com}\\$^{1}$Surya Sen Mahavidyalaya, Siliguri - 734004, West Bengal, India\\$^{2}$Department of Mathematics, University of North Bengal, \\Siliguri - 734013, West Bengal, India.}
\date{}
\maketitle

\begin{abstract}
A comparative study of the Homotopy Analysis method and an improved
Renormalization Group method is presented in the context of the Rayleigh and
the Van der Pol equations. Efficient approximate formulae as functions of the
nonlinearity parameter $\varepsilon$ for the amplitudes $a(\varepsilon)$ of
the limit cycles for both these oscillators are derived. The improvement in
the Renormalization group analysis is achieved by invoking the idea of
nonlinear time that should have significance in a nonlinear system. Good
approximate plots of limit cycles of the concerned oscillators are also
presented within this framework.

\end{abstract}

\begin{center}
\end{center}

\section{Introduction}

The study of non-perturbative methods for nonlinear differential equations is
of considerable recent interest. Among the various well known singular
perturbation techniques such as multiple scale analysis, method of boundary
layers, WKB method and so on \cite{Jordan Smith, Bender Orszag}, the recently
developed homotopy analysis method (HAM) \cite{Liao HAM, Lopez VDP Amp} and
the Renormalization group method (RGM) \cite{Chen G O RG 1994, Chen G O RG,
DeVille RG, Sarkar Bhattacharjee} appear to be very attractive. The aim of
these new improved methods is to derive in an unified manner uniformly valid
asymptotic quantities of interest for a given nonlinear dynamical problem.
Although formulated almost parallely over the past decades or so, relative
strength and weakness of these two approaches have yet to be investigated
systematically. The purpose of this paper is to undertake a comparative study
of HAM and RGM in the context of the Rayleigh equation%
\begin{equation}
\ddot{y}+\varepsilon\left(  \frac{1}{3}\dot{y}^{3}-\dot{y}\right)  +y=0
\label{Rayleigh}%
\end{equation}
and the Van der Pol equation%
\begin{equation}
\ddot{x}+\varepsilon~\dot{x}\left(  x^{2}-1\right)  +x=0 \label{VdP}%
\end{equation}
where the dots are used to designate the derivative with respect to time $t$.
The Rayleigh and the Van der Pol (VdP) equations represent two closely related
nonlinear systems and have found significant applications in the study of self
excited oscillations arising in biology, acoustics, robotics, engineering etc.
\cite{Acoustics, Robotics}. It is easy to observe that differentiating
$\left(  \text{\ref{Rayleigh}}\right)  $ with respect to time $t$ and putting
$\dot{y}\left(  t\right)  =x\left(  t\right)  $ we obtain $\left(
\text{\ref{VdP}}\right)  $. Both these systems have unique isolated periodic
orbit (limit cycle). The amplitude of a periodic oscillation $y(t)$ $\left(
\text{or }x\left(  t\right)  \right)  $ is generally defined by $\max
\left\vert {y}\left(  {t}\right)  \right\vert $ $\left(  \text{or }%
\max\left\vert {x}\left(  {t}\right)  \right\vert \right)  $ over the entire
cycle. It is well known that the naive perturbative solutions of these
equations are useful when $0<\varepsilon\ll1$ and yields the asymptotic value
$a(\varepsilon)\approx2$ of the amplitude for the limit cycle correctly. For
$\varepsilon\gg1$, simple analysis based on singular perturbation theory also
yields the asymptotic amplitude for the relaxation oscillation as
$a(\varepsilon)\approx2$ for the VdP equation. However, the conventional
perturbative approaches fail when $\varepsilon$ is finite. One of the aim of
this paper is to determine efficient approximate formulae for the amplitude of
the limit cycle for the above systems by both HAM and RGM. Lopez et al
\cite{Lopez VDP Amp} have reported an efficient formula for the amplitude of
the VdP limit cycle by HAM. We note here that a key difference in Rayleigh and
VdP oscillators is the fact that with increase in input energy (voltage), the
amplitude of the Rayleigh periodic oscillation increases, when that of the VdP
oscillator remains almost constant at the value 2, with possible increase in
the corresponding frequency only. For large $\varepsilon\ (\geq1)$ relaxation
oscillations, on the other hand, the Rayleigh system shows up a rather fast
building up and slow subsequent release of internal energy, when the VdP
models the reverse behaviour, with slow rise and fast drop in the accumulated energy.

As remarked above, HAM and RGM are formulated to determine the uniformly valid
global asymptotic behaviours of relevant dynamical quantities like amplitude,
period, frequency etc. related to periodic solutions of these equations for
finite values of $\varepsilon$, by devising efficient methods in eliminating
divergent secular terms of the naive perturbation theory. HAM seems to have
the advantage of yielding uniformly convergent solutions of very high order in
the nonlinearity parameter $\varepsilon$ utilizing a freedom in the choice of
a free parameter $h$. The computation of higher order term could be
facilitated by symbolic computational algorithms. This method is used to
obtain good approximate solutions for the VdP equation by a number of authors
\cite{Lopez VDP Amp, Liao VDP}. Lopez et al \cite{Lopez VDP Amp} derived
efficient formulae for estimating the amplitude of the limit cycle of the VdP
equation for all values of $\varepsilon>0$. Although, HAM is now considered to
be an efficient method in the study of non-perturbative asymptotic analysis,
it is recently pointed out \cite{Meijer} that this method might fail even in
some innocent looking nonlinear problems.

The RGM, on the other hand, has a rich history, being originally formulated
for managing divergences in the quantum field theory and later having deep
applications in phase transitions and critical phenomena in statistical
mechanics. Subsequently, Chen et al \cite{Chen G O RG 1994, Chen G O RG}
successfully translated the RG formalism into the study of nonlinear
differential equations. It is noted that RGM is more efficient and accurate
than conventional singular perturbative approaches in obtaining global
informations from a naive perturbation series in $\varepsilon$. It is also
recognized that RGM generated expansions yield $\varepsilon$-dependent
space/time scales naturally, when conventional approaches normally require
invoking such scales in an ad hoc manner. The pertubative RGM, however,
appears to have the limitation that the computation of higher order
renormalized solutions could be quite involved and tedious. More serious is
the inability of assuring the convergence of the renormalized expansions for
large nonlinearity parameter. Further, there is still no evidence in the
literature that RGM could be employed successfully to asymptotic estimation of
the amplitude of an isolated periodic orbit for all values of $\varepsilon$ as
was reported for HAM \cite{Lopez VDP Amp}.

Here we report analytic expressions of the amplitude of the periodic solutions
of both the Rayleigh equation $\left(  \text{\ref{Rayleigh}}\right)  $ and the
VdP equation $\left(  \text{\ref{VdP}}\right)  $ as functions of $\varepsilon
$. We have made a comparative study of these two sets of formulae using both
HAM and RGM. The HAM contains a control parameter $h=h\left(  \varepsilon
\right)  $ which controls the convergence of the approximation to the
numerically computed exact value of the amplitude for all values of
$\varepsilon$. Suitable choice of $h$ can control the relative percentage error.

The original RGM gives an approximation to the exact solution for small values
of $\varepsilon$. We report here the RG solution upto order $3$. To the
authors' knowledge this seems to be the first higher order computation other
than second order computations reported so far by various authors \cite{Chen G
O RG, Sarkar Bhattacharjee}. A comparison of the amplitude of the periodic
cycle with the exact computations reveals that even the present higher order
perturbative approximations fails to give accurate estimation for moderate
values of $\varepsilon$. As the higher order RG computations are quite
laborious and inefficient, it is very unlikely that higher order computations
of amplitude would improve the quality of the estimated amplitude of the limit
cycle. Further, in RGM one does not have the resource of a free parameter
equivalent to $h(\varepsilon)$ of HAM to improve the convergence of the RG expansions.

A major contribution in the present study is to propose an \emph{improved} RGM
(IRGM). In IRGM, we advocate the concept of \emph{nonlinear time} \cite{DPD
Raut, DPD, DPD Sen} that extends the original RG idea of eliminating the
divergent secular term of the form $(t-t_{0})\sin t$, where $t_{0}$ is the
initial time, of the naive perturbation series for the solution of the
nonlinear problem, by exploiting the arbitrariness in fixing the initial
moment $t_{0}$. The original prescription rests on introducing new initial
time $\tau$ in the form $(t-\tau+\tau-t_{0})\sin t$ and to allow the
renormalized amplitude $R=R(\tau)$ and phase $\theta=\theta(\tau)$ of the
renormalized solution to depend on the new parameter, viz., $\tau-t_{0}$ so
that the original naive perturbative, constant values of amplitude $R_{0}$ and
phase $\theta_{0}\ (=0)$ (say) `flow' following the RG flow equations of the
form
\begin{equation}
\frac{dR}{d\tau}=f\left(  R,\varepsilon\right)  ,\ \frac{d\theta}{d\tau
}=g\left(  R,\varepsilon\right)  \label{rgflow}%
\end{equation}
The functions in the right hand sides of the RG flow equations, in general,
should depend both on $R$ and $\theta$, besides the explicit $\varepsilon$
dependence. We suppress the $\theta$ dependence for simplicity that should
suffice for our present analysis of the Rayleigh and VdP equations $($c.f.
equations $\left(  \text{\ref{Amp Eq Order 3}}\right)  $, $\left(
\text{\ref{Phase Eq Order 3}}\right)  )$. The flow equations are derived from
the consistency condition that the actual renormalized solution $y\left(
t,\tau\right)  $ should be independent of the arbitrary initial adjustment
$\tau$: $\frac{\partial y}{\partial\tau}=0$. The final form of the uniformly
valid RG solution $y_{R}\left(  t\right)  $ is obtained by setting $\tau=t$
that eliminates the secular terms. Let us remark here that the actual
convergence of the RG expansions is not well addressed and should require
further investigations. Moreover, estimation of asymptotic amplitude for a
limit cycle as $t\rightarrow\infty$, for instance, from the perturbation
expansion of $f$ is expected to fail for $\varepsilon>\approx O\left(
1\right)  $.

In the framework of nonlinear time, we suppose the arbitrary initial time
$\tau$ to depend explicitly on the nonlinearity parameter (coupling strength)
$\varepsilon$ of the nonlinear equation, so that one can write $\tau
/\varepsilon=\varepsilon^{h}$ where $h=h\left(  \varepsilon t\right)
,\ \varepsilon t>1$ is a slowly varying (almost constant), free (asymptotic)
control parameter for $t\rightarrow\infty$ and $\varepsilon\ \rightarrow$
either to 0 or $\infty$, to be utilized judiciously to improve the convergence
and non-perturbative global asymptotic behaviour of the original RG proposal
($h<0$ for $0<\varepsilon<1$). In Appendix, we give an overview, in brief, of
an extended analytic framework that naturally supports nontrivial existence of
such an asymptotic scaling parameter $h(\tilde{\tau})$ as a function of the
rescaled $O(1)$ variable $\tilde{\tau}=\varepsilon t\sim O(1)$, satisfying
what we call the \emph{principle of duality structure}. The secular terms in
the naive perturbation series would now be altered instead as $\left(
t-\tau/\varepsilon+\tau/\varepsilon-t_{0}\right)  \sin t$ and we obtain the
new RG flow equations in the form
\begin{equation}
\frac{dR}{d\tau}=f_{0}\left(  R\right)  \left(  1+O\left(  \varepsilon
^{2}\right)  \right)  ,\ \frac{d\theta}{d\tau}=\varepsilon g_{1}\left(
R\right)  \left(  1+O\left(  \varepsilon^{3}\right)  \right)  \label{rgflown}%
\end{equation}
where $f_{0}(R)$ and $g_{1}(R)$ are nonzero, minimal order $R$ dependent terms
in the respective perturbation series. Following the analogy of RG
prescription in annulling secular divergence through corresponding `flowing'
of the renormalized perturbative amplitude and phase, we next make \emph{the
key assumption that there exists, for a given nonlinear oscillation, a set of
right control parameters $h^{i}$ that would absorb any possible secular or
other kind of divergence in the higher order perturbation series}, (see
Appendix for justification), so that in the asymptotic limit $t\rightarrow
\infty$, one obtains the finite, non-perturbative flow equations directly for
the periodic orbit of the nonlinear system
\begin{equation}
\frac{da}{d\tau_{1}}=f_{0}(a),\ \frac{d\theta}{d\tau_{2}}=g_{1}%
(a)\label{rgflownp}%
\end{equation}
where $\tau_{i}=\varepsilon^{ih_{RG}^{i}(\tilde{\tau})}$, and $h_{RG}%
^{i}(\tilde{\tau})$ is a finite \emph{scale} independent control parameter in
the rescaled variable $\tilde{\tau}\sim$ $O(1)$ and $a(\varepsilon
)=\underset{t\rightarrow\infty}{\lim}R(\varepsilon,t)$ is the $\varepsilon$-
dependent amplitude of the limit cycle. A simple quadrature formula should
then relate the control parameter $h_{RG}=h_{RG}^{1}$ with the amplitude
$a(\varepsilon)$. As a consequence, adjusting the control parameter $h_{RG}$
suitably, one can generate an efficient algorithm to estimate the amplitude
$a(\varepsilon)$ that would compare well with the exact values, upto any
desired accuracy. It will transpire that the control parameter $h_{RG}%
(\varepsilon)$ must respect some asymptotic conditions depending on the
characteristic features of a particular relaxation oscillation $($c.f. Section
\ref{Sec M-RG Sol Description}$)$.

Exploiting the rescaling symmetry, one may as well rewrite the above
non-perturbative flow equations (\ref{rgflownp}) in the equivalent
$\tilde{\tau}\sim O(1)$-dependent scaling variable $\tau={\tilde{\tau}%
}^{H_{RG}(\tilde{\tau})}$, $\left(  H_{RG}(\tilde{\tau})=h_{RG}(\tilde{\tau
})\frac{\log\tilde{\tau}}{\log\varepsilon}\right)  $, for each fixed value of
the nonlinearity parameter $\varepsilon$ that should expose small scale
$\tilde{\tau}\sim O(1)$-dependent variation of the amplitude. As a biproduct
that would allow one to retrieve an efficient approximation of the limit cycle
orbit for the nonlinear oscillator. It turns out that the general framework of
IRGM is quite successful in obtaining excellent fits for the limit cycle orbit
even for relaxation oscillation corresponding to nonlinearity parameters
$\varepsilon\geq1$.

It follows that the application of the of idea of nonlinear time 
in RG formalism offers one with a robust formalism
for global asymptotic analysis for a general nonlinear system that might even
be advantageous in many respects compared to HAM. The application of nonlinear
time in HAM will be considered separately.

The paper is organized as follows. In Section \ref{Sec HAM Sol} we have
deduced the solution to the equation $\left(  \text{\ref{Rayleigh}}\right)  $
by HAM. In Section \ref{Sec RG Sol} we compute the classical RG solution upto
$O(\varepsilon^{3})$ order and compare estimated values of the limit cycle
amplitude with the exact values. The improved RG method is presented in
Section \ref{Sec M-RG Sol Description}. This introduces a control parameter
$h_{RG}$ in the RG analysis. In Subsection \ref{SubSec Amplitude Estimation}
approximate analytic formulae are deduced for the amplitudes of the limit
cycles of the Rayleigh and the VdP equations. Efficient match with the exact
values can be obtained by appropriate choice of $h_{RG}$. In Subsection 4.2 we
present the efficient of approximate limit cycle orbits for the Rayleigh and VdP
oscillators for $\varepsilon=5$. We close our discussions in Sec. 5. In Appendix 1,
 we present a brief outline of the formal structure of the analytic formalism presented here. 
In Appendix 2, an alternative approach in the derivation of non-perturbative flow equations is presented.

\section{Computation of Amplitude by HAM\label{Sec HAM Sol}}

The Homotopy Analysis method proposed by Liao \cite{Liao HAM, Liao VDP} is
used to obtain the solution of non-linear equation even if the problem does
not contain a small or large parameter. HAM always gives a family of functions
at any given order of approximation. An auxiliary parameter $h$ is introduced
in HAM to control the convergence region of approximating series involved in
this method to the exact solution. HAM is based on the idea of homotopy in
topology. In simple language, it involves continuous deformation of the
solution of a linear ordinary differential equation (ODE) to that of desired
nonlinear ODE. The solution of linear ODE gives a set of functions called
\textit{base functions}. One advantage of HAM is that it can be used to
approximate a nonlinear problem by efficient choice of different sets of base
functions. A suitable choice of the set of base functions and the convergence
control parameter can speed up the convergence process.

In this paper we consider the self-excited system $\left(
\text{\ref{Rayleigh}}\right)  $, which can be written as the ODE%
\begin{equation}
\ddot{U}\left(  t\right)  +\varepsilon\left(  \frac{1}{3}\dot{U}^{3}\left(
t\right)  -\dot{U}\left(  t\right)  \right)  +U\left(  t\right)  =0,\quad
t\geq0 \label{RL Eq}%
\end{equation}
where the dot denotes the derivative with respect to the time $t$. A limit
cycle represents an isolated periodic motion of a self-excited system. This is
an isolated closed curve $\Gamma$ $\left(  \text{say}\right)  $ in the phase
plane so that any path in its suitable small neighbourhood starting from a
point, specified by some given initial condition, ultimately converges to
$($or diverge from$)$ $\Gamma$. Consequently, this periodic motion represented
by limit cycle is independent of initial conditions. It, however, involves the
frequency $\omega$ and the amplitude $a$ of the oscillation. Therefore,
without loss of generality, we consider an initial condition%
\begin{equation}
U\left(  0\right)  =a,\qquad\dot{U}\left(  0\right)  =0\text{.} \label{IC}%
\end{equation}

In \cite{Chen G O RG}, an alternative initial condition i.e. $U(0)=0,\ \dot
{U}(0)=a$ was considered. Let, with slight abuse of notations,
\[
\tau=\omega t\text{ and }U\left(  t\right)  =a~u\left(  \tau\right)  \text{.}%
\]
so that $\left(  \text{\ref{RL Eq}}\right)  $ and $\left(  \text{\ref{IC}%
}\right)  $ respectively become%
\begin{equation}
\omega^{2}u^{\prime\prime}\left(  \tau\right)  +\varepsilon\left(  \frac{1}%
{3}a^{2}\omega^{2}u^{\prime2}\left(  \tau\right)  -1\right)  \omega u^{\prime
}\left(  \tau\right)  +u\left(  \tau\right)  =0 \label{Normalized RL}%
\end{equation}
and%
\begin{equation}
u\left(  0\right)  =1\text{,\quad}u^{\prime}\left(  0\right)  =0\text{.}
\label{Normalized-IC}%
\end{equation}
Since the limit cycle represents a periodic motion, so we suppose that the
initial approximation to the solution $u\left(  \tau\right)  $ to $\left(
\text{\ref{Normalized RL}}\right)  $ can be taken as%
\[
u_{0}\left(  \tau\right)  =\cos\tau
\]
Let, $\omega_{0}$ and $a_{0}$ respectively denote the initial approximations
of the frequency $\omega$ and the amplitude $a$.

We consider a linear operator%
\begin{equation}
\mathcal{L}\left[  \phi\left(  \tau,p\right)  \right]  =\omega_{0}^{2}\left[
\frac{\partial^{2}\phi\left(  \tau,p\right)  }{\partial\tau^{2}}+\phi\left(
\tau,p\right)  \right]  \label{L Operator}%
\end{equation}
so that for the coefficients $C_{1}$ and $C_{2}$%
\begin{equation}
\mathcal{L}\left(  C_{1}\sin\tau+C_{2}\cos\tau\right)  =0 \label{L IC}%
\end{equation}
We further consider a nonlinear operator%
\begin{align}
&  \mathcal{N}\left[  \phi\left(  \tau,p\right)  ,\Omega\left(  p\right)
,A\left(  p\right)  \right] \nonumber\\
&  =\Omega^{2}\left(  p\right)  \frac{\partial^{2}\phi\left(  \tau,p\right)
}{\partial\tau^{2}}+\varepsilon\left[  \frac{1}{3}A^{2}\left(  p\right)
\Omega^{3}\left(  p\right)  \left(  \frac{\partial\phi\left(  \tau,p\right)
}{\partial\tau}\right)  ^{3}-\Omega\left(  p\right)  \left(  \frac
{\partial\phi\left(  \tau,p\right)  }{\partial\tau}\right)  \right]
+\phi\left(  \tau,p\right)  \text{.} \label{NL Operator}%
\end{align}
Next, we construct a homotopy as%
\begin{equation}
\mathcal{H}\left[  \phi\left(  \tau,p\right)  ,h,p\right]  =\left(
1-p\right)  \mathcal{L}\left[  \phi\left(  \tau,p\right)  -u_{0}\left(
\tau\right)  \right]  -h~p\mathcal{~N}\left[  \phi\left(  \tau,p\right)
,\Omega\left(  p\right)  ,A\left(  p\right)  \right]  \label{Homotopy}%
\end{equation}
where $p\in\left[  0,1\right]  $ is the embedding parameter and $h$ a non-zero
auxiliary (control) parameter used to improve the convergence of series
expansions. Setting $\mathcal{H}\left[  \phi\left(  \tau,p\right)
,h,p\right]  =0$ we obtain zero-th order deformation equation%
\begin{equation}
\left(  1-p\right)  \mathcal{L}\left[  \phi\left(  \tau,p\right)
-u_{0}\left(  \tau\right)  \right]  -h~p\mathcal{~N}\left[  \phi\left(
\tau,p\right)  ,\Omega\left(  p\right)  ,A\left(  p\right)  \right]  =0
\label{Zero Deform Eq}%
\end{equation}
subject to the initial conditions%
\begin{equation}
\phi\left(  0,p\right)  =1,\quad\left.  \frac{\partial\phi\left(
\tau,p\right)  }{\partial\tau}\right\vert _{\tau=0}=0\text{.}
\label{Zero Deform IC}%
\end{equation}
Clearly, as $p$ increases from $p=0$ to $p=1$, $\left(
\text{\ref{Zero Deform Eq}}\right)  $ changes from $\mathcal{L}\left[
\phi\left(  \tau,p\right)  -u_{0}\left(  \tau\right)  \right]  =0$ to
$\mathcal{N}\left[  \phi\left(  \tau,p\right)  ,\Omega\left(  p\right)
,A\left(  p\right)  \right]  =0$ and as a consequence $\phi\left(
\tau,p\right)  $ varies from the initial guess $\phi\left(  \tau,0\right)
=u_{0}\left(  \tau\right)  =\cos\tau$ to the exact solution $\phi\left(
\tau,1\right)  =u\left(  \tau\right)  $, so does $\Omega\left(  p\right)  $
from $\omega_{0}$ to exact frequency $\omega$ and $A\left(  p\right)  $ from
$a_{0}$ to the exact amplitude $a$. It can be shown that assuming $\phi\left(
\tau,p\right)  $, $\Omega\left(  p\right)  $, $A\left(  p\right)  $ analytic
in $p\in\left[  0,1\right]  $ so that%
\begin{equation}
u_{k}\left(  \tau\right)  =\frac{1}{k!}\left.  \frac{\partial^{k}}{\partial
p^{k}}\phi\left(  \tau,p\right)  \right\vert _{p=0},\quad\omega_{k}=\frac
{1}{k!}\left.  \frac{\partial^{k}}{\partial p^{k}}\Omega\left(  p\right)
\right\vert _{p=0},\quad a_{k}=\frac{1}{k!}\left.  \frac{\partial^{k}%
}{\partial p^{k}}A\left(  p\right)  \right\vert _{p=0}
\label{k Deform Derivative}%
\end{equation}
we have,%
\begin{align}
u\left(  \tau\right)   &  =%
{\textstyle\sum\limits_{k=0}^{\infty}}
u_{k}\left(  \tau\right) \label{Sol Series}\\
\omega &  =%
{\textstyle\sum\limits_{k=0}^{\infty}}
\omega_{k}\label{Freq Series}\\
a  &  =%
{\textstyle\sum\limits_{k=0}^{\infty}}
a_{k} \label{Amp Series}%
\end{align}
where $u_{k}\left(  \tau\right)  $ are solutions of the $k$-th order
deformation equation%
\begin{equation}
\mathcal{L}\left[  u_{k}\left(  \tau\right)  -\chi_{k}u_{k-1}\left(
\tau\right)  \right]  =h~R_{k}\left(  \tau\right)  \label{k Deform Eq}%
\end{equation}
subject to the initial conditions%
\begin{equation}
u_{k}\left(  0\right)  =0,\qquad u_{k}^{\prime}\left(  0\right)  =0
\label{k Deform IC}%
\end{equation}
in which%
\begin{align}
R_{k}\left(  \tau\right)   &  =\frac{1}{\left(  k-1\right)  !}\left.
\frac{\partial^{k-1}}{\partial p^{k-1}}\mathcal{N}\left[  \phi\left(
\tau,p\right)  ,\Omega\left(  p\right)  ,A\left(  p\right)  \right]
\right\vert _{p=0}\nonumber\\
&  =%
{\textstyle\sum\limits_{n=0}^{k-1}}
u_{k-1-n}^{\prime\prime}\left(  \tau\right)
{\textstyle\sum\limits_{j=0}^{n}}
\omega_{j}\omega_{n-j}+u_{k-1}\left(  \tau\right) \nonumber\\
&  +\frac{\varepsilon}{3}%
{\textstyle\sum\limits_{n=0}^{k-1}}
{\textstyle\sum\limits_{i=0}^{n}}
\left(
{\textstyle\sum\limits_{r=0}^{i}}
a_{r}a_{i-r}\right)  \times\left(
{\textstyle\sum\limits_{s=0}^{n-i}}
\omega_{s}%
{\textstyle\sum\limits_{h=0}^{n-i-s}}
\omega_{h}\omega_{n-i-s-h}\right) \nonumber\\
&  \times\left(
{\textstyle\sum\limits_{j=0}^{k-1-n}}
u_{j}^{\prime}\left(  \tau\right)
{\textstyle\sum\limits_{m=0}^{k-1-n-j}}
u_{m}^{\prime}\left(  \tau\right)  u_{k-1-n-j-m}^{\prime}\left(  \tau\right)
\right)  -\varepsilon%
{\textstyle\sum\limits_{n=0}^{k-1}}
\omega_{n}u_{k-1-n}^{\prime}\left(  \tau\right)  \label{k Deform RHS}%
\end{align}
and%
\begin{equation}
\chi_{k}=\left\{
\begin{array}
[c]{cc}%
0, & k\leq1,\\
1, & k>1.
\end{array}
\right.  \label{Chi}%
\end{equation}
To ensure that the solution to the $k$-th order deformation equation $\left(
\text{\ref{k Deform Eq}}\right)  $ do not contain the secular terms $\tau
\sin\tau$ and $\tau\cos\tau$ the coefficients of $\sin\tau$ and $\cos\tau$ in
the expressions of $R_{k}$ in $\left(  \text{\ref{k Deform RHS}}\right)  $
must vanish\ giving successive values of $\omega_{k}$ and $a_{k}$.

The linear equation $\mathcal{L}\left(  \phi\left(  \tau,p\right)  \right)
=0$ represents a simple harmonic motion with frequency $1$. So, we choose the
initial guess of $\omega$ as $\omega_{0}=1$. Again, by perturbation method
\cite{Jordan Smith} we find $a\rightarrow2$ as $\varepsilon\rightarrow0$. So,
we choose the initial guess of $a$ as $a_{0}=2$. Solving the differential
equations given by $\left(  \text{\ref{Zero Deform Eq}}\right)  $, $\left(
\text{\ref{Zero Deform IC}}\right)  $, $\left(  \text{\ref{k Deform Eq}%
}\right)  $, $\left(  \text{\ref{k Deform IC}}\right)  $ and avoiding the
generation of secular terms in each iteration we obtain%
\[
u_{1}\left(  \tau\right)  =-\frac{1}{24}h\varepsilon\sin3\tau+\frac{1}%
{8}h\varepsilon\sin\tau\text{, }\omega_{1}=-\frac{1}{16}h\varepsilon
^{2}\text{, }a_{1}=\frac{1}{8}h\varepsilon^{2}%
\]%
\begin{align*}
u_{2}\left(  \tau\right)   &  =\left(  \frac{1}{384}h^{2}\varepsilon^{3}%
-\frac{1}{24}h^{2}\varepsilon-\frac{1}{24}h\varepsilon\right)  \sin3\tau
-\frac{1}{64}h^{2}\varepsilon^{2}\cos3\tau\\
&  +\frac{1}{64}h^{2}\varepsilon^{2}\cos\tau+\left(  \frac{1}{8}%
h^{2}\varepsilon-\frac{1}{128}h^{2}\varepsilon^{3}+\frac{1}{8}h\varepsilon
\right)  \sin\tau
\end{align*}
so that%
\begin{align*}
R_{1}  &  =\frac{1}{3}\varepsilon\sin3\tau\\
R_{2}  &  =\frac{1}{24}\left[  3h\varepsilon^{2}\cos3\tau+\left(
8h\varepsilon-\frac{1}{2}h\varepsilon^{3}\right)  \sin3\tau\right]
\end{align*}
Computing $R_{k}$ successively, we can find the successive expressions of
$u_{k}\left(  \tau\right)  $, $\omega_{k}$ and $a_{k}$. The first order
approximation to the amplitude in $\left(  \text{\ref{Amp Series}}\right)  $
is%
\begin{equation}
a\approx a_{0}+a_{1}=2+\frac{1}{8}h\varepsilon^{2}=a_{E}\left(  \varepsilon
\right)  \text{ }\left(  \text{say}\right)  \text{.}
\label{Amp Approx 1 Order}%
\end{equation}
\newline\begin{figure}[ptb]
\begin{center}
\includegraphics[width=3in]{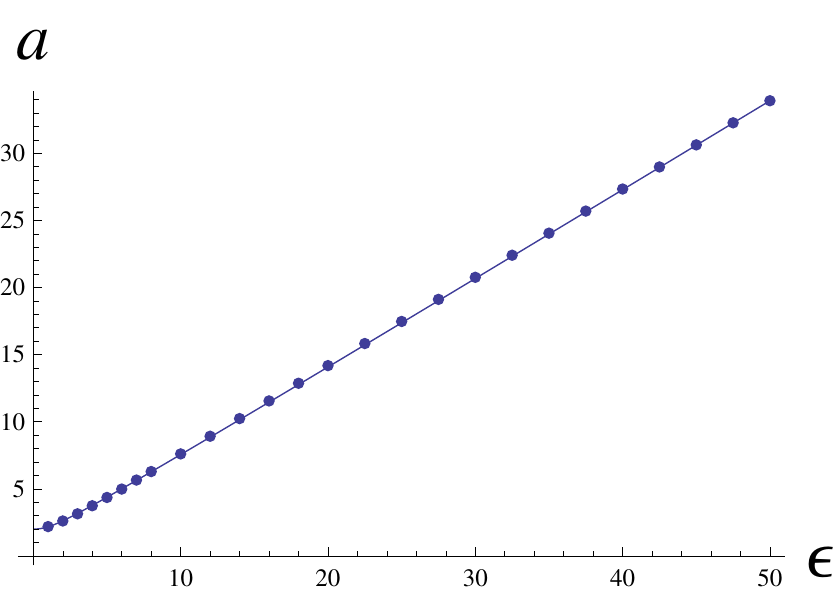}
\end{center}
\caption{The exact amplitude of Rayleigh Equation $($by solid line$)$ and its
approximation $a_{E}\left(  \varepsilon\right)  $ given by $\left(
\text{\ref{HAM Amp New}}\right)  $ $($by bold points$)$ for $0<\varepsilon
\leq50$.}%
\label{Fig Exact Amp RL}%
\end{figure}

The above first order expression for the amplitude involves as yet arbitrary
control parameter $h$. Lopez et al \cite{Lopez VDP Amp} proposed specific
$\varepsilon$-dependent expressions for $h$ to obtain an efficient formula for
the VdP limit cycle amplitude. They made the proposal that $h$, besides being
continuous, must also vanish in the limits of $\varepsilon\rightarrow0$ and
$\varepsilon\rightarrow\infty$ to reproduce the zeroth order perturbative
solutions. In our application of HAM for the Rayleigh limit cycle amplitude,
we have chosen a different set of base functions and so can weaken the
condition considerably, both on the continuity and the asymptotic limit
$\varepsilon\rightarrow\infty$. From careful inspections of the graph of the
exact amplitude (Fig.1), it turns out that an appropriate ansatz for the
control parameter $h$ is given by
\begin{equation}
h=\frac{1}{0.5+\varepsilon~b\left(  \varepsilon\right)  }
\label{Control Parameter}%
\end{equation}
where, $b\left(  \varepsilon\right)  $ is taken as the step function in the
domain $0<\varepsilon\leq50$ as follows%
\[%
\begin{tabular}
[c]{rccccc}%
$\varepsilon:$ & $0<\varepsilon\leq4$ & $4<\varepsilon\leq5$ & $5<\varepsilon
\leq7$ & $7<\varepsilon\leq8$ & $8<\varepsilon\leq9$\\
$b\left(  \varepsilon\right)  :$ & $0.162$ & $0.165$ & $0.168$ & $0.171$ &
$0.174$\\
$\varepsilon:$ & $9<\varepsilon\leq11$ & $11<\varepsilon\leq15$ &
$15<\varepsilon\leq20$ & $20<\varepsilon\leq30$ & $30<\varepsilon\leq50$\\
$b\left(  \varepsilon\right)  :$ & $0.176$ & $0.179$ & $0.181$ & $0.183$ &
$0.185$%
\end{tabular}
\
\]
With this particular form of $h$, we are able to find an analytic
approximation $a_{E}\left(  \varepsilon\right)  $ to the numerically computed
exact value $a=a\left(  \varepsilon\right)  $ in the domain $0<\varepsilon
\leq50$ with maximum relative percentage error $\left\vert \dfrac{a_{E}\left(
\varepsilon\right)  -a\left(  \varepsilon\right)  }{a\left(  \varepsilon
\right)  }\times100\right\vert $ less than $1\%$. Obviously, better accuracy
fit can be obtained by considering finer subdivisions in the definition of
$b(\varepsilon)$. We remark that a piece-wise continuous $\varepsilon$
dependence of $h$ as above is admissible in the framework of HAM.

Since the exact graph of $a(\varepsilon)$ is almost a straight line for
sufficiently large $\varepsilon$ $\left(  7<\varepsilon\leq50\right)  $, we
can reduce the number of steps to 4 only. Let us choose%

\begin{equation}
h=\frac{8m}{\varepsilon}-\frac{56m}{\varepsilon^{2}}+\frac{8c}{\varepsilon
^{2}}-\frac{16}{\varepsilon^{2}},\qquad7<\varepsilon\leq50
\label{Control Parameter New}%
\end{equation}
so that $\left(  \text{\ref{Amp Approx 1 Order}}\right)  $ becomes%
\begin{equation}
a_{E}\left(  \varepsilon\right)  =\left\{
\begin{array}
[c]{lc}%
2+\frac{1}{8}\left(  \frac{1}{0.5+0.162~\varepsilon}\right)  \varepsilon^{2} &
0<\varepsilon\leq4\\
2+\frac{1}{8}\left(  \frac{1}{0.5+0.165~\varepsilon}\right)  \varepsilon^{2} &
4<\varepsilon\leq5\\
2+\frac{1}{8}\left(  \frac{1}{0.5+0.168~\varepsilon}\right)  \varepsilon^{2} &
5<\varepsilon\leq7\\
m\left(  \varepsilon-7\right)  +c\text{,} & 7<\varepsilon\leq50
\end{array}
\right.  \label{HAM Amp New}%
\end{equation}
where $m$ and $c$ are computed from the exact solution as%
\[
m=\dfrac{a\left(  50\right)  -a\left(  7\right)  }{50-7}=0.657692\text{ and
}c=a\left(  7\right)  =5.63108
\]
keeping the maximum relative percentage error $\left\vert \dfrac{a_{E}\left(
\varepsilon\right)  -a\left(  \varepsilon\right)  }{a\left(  \varepsilon
\right)  }\times100\right\vert $ less than $1\%$. The plot of $a_{E}\left(
\varepsilon\right)  $ given by $\left(  \text{\ref{HAM Amp New}}\right)  $ is
shown by bold points in Figure \ref{Fig Exact Amp RL} (explicit
discontinuities of $h$ at $\varepsilon=4,5$ and $7$ are not visible at the
resolution of the plotted figure) . \begin{figure}[h]
\begin{center}
\includegraphics[width=3in]{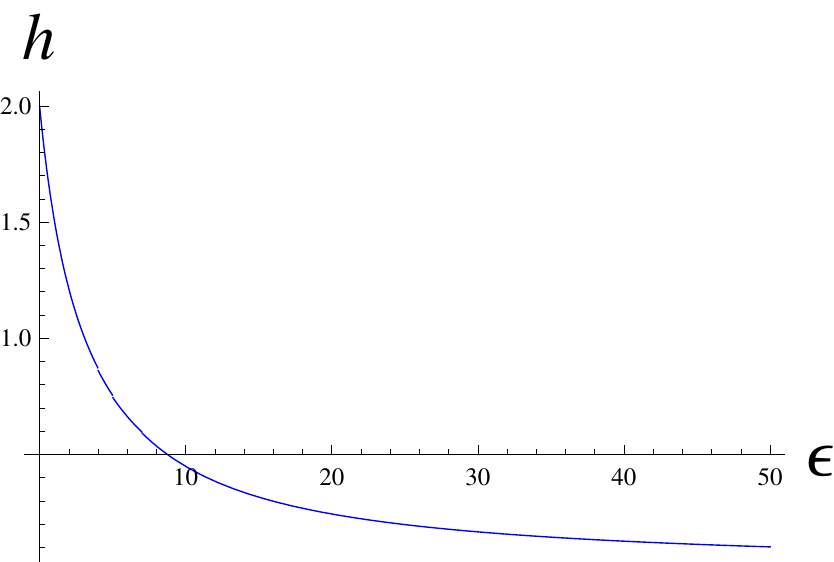}
\end{center}
\caption{The graph of $h\left(  \varepsilon\right)  $ used for approximation
of the amplitude by HAM given by $\left(  \text{\ref{HAM Amp New}}\right)  $
for $0<\varepsilon\leq50$.}%
\label{Fig h Ham Rayleigh(epsilon)}%
\end{figure}As remarked above, Lopez et. al. \cite{Lopez VDP Amp} proposed
that a reasonable property for $h$ would be to vanish in the limits as
$\varepsilon\rightarrow0$ and $\varepsilon\rightarrow\infty$. However, from
$\left(  \text{\ref{Control Parameter}}\right)  $ and $\left(
\text{\ref{Control Parameter New}}\right)  $ we observe that a suitable
approximation to the amplitude of Rayleigh equation can be obtained even if
$h$ do not follow this property. The graph of $h\left(  \varepsilon\right)  $
is given in Figure \ref{Fig h Ham Rayleigh(epsilon)} for $0<\varepsilon\leq50$
(discontinuity in $h$ is not visible at the level of resolution in the figure).

To summarize, one can obtain more accurate approximate formula by suitable
choices of the control parameter $h\left(  \varepsilon\right)  $ upto any
desired level of accuracy. We also note that a piecewise continuous control
parameter $h$ enables us to obtain good approximation by solving only the
first order deformation equation. However, the first order HAM estimated
amplitude $a(\varepsilon)$ is $O\left(  \varepsilon^{2}\right)  $.

We report the estimation of the amplitude of the limit cycle for the Rayleigh
and VdP equations by the improved RG method in Subsection
\ref{SubSec Amplitude Estimation}. We do not undertake the computation of the
VdP amplitude by HAM separately, as that was already reported by Lopez et al
\cite{Lopez VDP Amp}.

\section{Computation of Amplitude by RG Method\label{Sec RG Sol}}

The Renormalization Group method (RGM) introduced by Chen, Goldenfeld and Oono
(CGO) \cite{Chen G O RG 1994, Chen G O RG} gives a unified formal approach to
derive asymptotic expansions for the solutions of a large class of nonlinear
ODEs. The RG method is used in solid state physics, quantum field theory and
other areas of physics. One advantage of RGM is that it starts from naive
perturbation expansion of a problem and is expected to yield automatically the
gauge functions such as fractional powers of $\varepsilon$ and logarithmic
terms in $\varepsilon$ in the renormalized expansion. One does not require to
have any prior knowledge to prescribe these unexpected gauge functions in an
ad hoc manner. DeVille et. al. \cite{DeVille RG}\ have introduced an
algorithmic approach for RGM which we adopt for the following application. As
it will transpire the RGM appears to be deficient in estimating amplitude of a
periodic orbit because of the absence of any free control parameter. In a
latter section we have improved this RGM to incorporate a control parameter
similar to HAM and derive efficient estimations of amplitudes of both the
Rayleigh and VdP equations. However, before the introduction of the improved
RG method (IRGM), we first discuss the \textit{conventional} RG method, given
by DeVille et. al. and use it to obtain amplitude and phase equations for the
Rayleigh equation $\left(  \text{\ref{Rayleigh}}\right)  $. These equations
are already obtained in \cite{Chen G O RG, DeVille RG} to the order $O\left(
\varepsilon^{3}\right)  $ which agree with the experimental values as
$\varepsilon\rightarrow0$ only. We have extended these results to the order
$O\left(  \varepsilon^{4}\right)  $ and notice that higher order perturbative
computations of the RG flow equations would fail to obtain good estimation of
the amplitude of the periodic cycle for all values of $\varepsilon$.

Substituting the naive expansion%
\[
y\left(  t\right)  =y_{0}\left(  t\right)  +\varepsilon y_{1}\left(  t\right)
+\varepsilon^{2}y_{2}\left(  t\right)  +\varepsilon^{3}y_{3}\left(  t\right)
+\cdots
\]
in $\left(  \text{\ref{Rayleigh}}\right)  $, we find at each order%
\begin{align*}
O\left(  1\right)   &  :\ddot{y}_{0}+y_{0}=0\\
O\left(  \varepsilon\right)   &  :\ddot{y}_{1}+y_{1}=\dot{y}_{0}-\frac{1}%
{3}\dot{y}_{0}^{3}\\
O\left(  \varepsilon^{2}\right)   &  :\ddot{y}_{2}+y_{2}=\dot{y}_{1}-\dot
{y}_{0}^{2}\dot{y}_{1}\\
O\left(  \varepsilon^{3}\right)   &  :\ddot{y}_{3}+y_{3}=\dot{y}_{2}-\dot
{y}_{0}^{2}\dot{y}_{2}-\dot{y}_{1}^{2}\dot{y}_{0}%
\end{align*}
The solutions are\newline$y_{0}\left(  t\right)  =Ae^{i\left(  t-t_{0}\right)
}+c.c.\bigskip$\newline$y_{1}\left(  t\right)  =\frac{1}{24}iA^{3}e^{i\left(
t-t_{0}\right)  }+\frac{1}{2}A\left(  1-AA^{\ast}\right)  \left(
t-t_{0}\right)  e^{i\left(  t-t_{0}\right)  }-\frac{1}{24}iA^{3}e^{3i\left(
t-t_{0}\right)  }+c.c\bigskip$\newline$y_{2}\left(  t\right)  =\left(
\frac{1}{32}A^{3}-\frac{3}{64}A^{4}A^{\ast}\right)  e^{i\left(  t-t_{0}%
\right)  }\bigskip$\newline$\left.  {}\right.  \hspace{0.5in}+\left(
-\frac{1}{24}iA^{4}A^{\ast}+\frac{1}{16}iA^{3}(A^{\ast})^{2}+\frac{1}%
{48}iA^{3}+\frac{1}{48}iA^{2}(A^{\ast})^{3}-\frac{1}{8}iA\right)  \left(
t-t_{0}\right)  e^{i\left(  t-t_{0}\right)  }\bigskip$\newline$\left.
{}\right.  \hspace{0.5in}+\left(  \frac{3}{8}A^{3}(A^{\ast})^{2}-\frac{1}%
{2}A^{2}A^{\ast}+\frac{1}{8}A\right)  \left(  t-t_{0}\right)  ^{2}e^{i\left(
t-t_{0}\right)  }+\left(  \frac{3}{64}A^{4}A^{\ast}-\frac{1}{32}A^{3}+\frac
{1}{192}A^{5}\right)  e^{3i\left(  t-t_{0}\right)  }\bigskip$\newline$\left.
{}\right.  \hspace{0.5in}-\frac{1}{16}iA^{3}\left(  1-AA^{\ast}\right)
\left(  t-t_{0}\right)  e^{3i\left(  t-t_{0}\right)  }-\frac{1}{192}%
A^{5}e^{5i\left(  t-t_{0}\right)  }+c.c.\bigskip$\newline$y_{3}\left(
t\right)  =\left(  -\frac{1}{384}iA^{6}A^{\ast}+\frac{37}{1536}iA^{5}(A^{\ast
})^{2}+\frac{1}{2304}iA^{5}+\frac{1}{512}iA^{4}(A^{\ast})^{3}-\frac{7}%
{256}iA^{4}A^{\ast}-\frac{1}{128}iA^{3}\right)  e^{i\left(  t-t_{0}\right)
}\bigskip$\newline$\left.  {}\right.  \hspace{0.5in}+\left(
\begin{array}
[c]{c}%
+\frac{1}{1152}A^{6}A^{\ast}+\frac{5}{128}A^{5}(A^{\ast})^{2}-\frac{119}%
{1152}A^{4}(A^{\ast})^{3}+\frac{11}{384}A^{3}(A^{\ast})^{4}\bigskip\\
-\frac{7}{128}A^{4}A^{\ast}+\frac{1}{48}A^{3}+\frac{11}{64}A^{3}(A^{\ast}%
)^{2}-\frac{1}{64}A^{2}(A^{\ast})^{3}%
\end{array}
\right)  \left(  t-t_{0}\right)  e^{i\left(  t-t_{0}\right)  }\bigskip
$\newline$\left.  {}\right.  \hspace{0.5in}+\left(
\begin{array}
[c]{c}%
+\frac{3}{64}iA^{5}(A^{\ast})^{2}-\frac{3}{32}iA^{4}(A^{\ast})^{3}-\frac
{1}{24}iA^{4}A^{\ast}-\frac{1}{32}iA^{3}(A^{\ast})^{4}\bigskip\\
+\frac{3}{32}iA^{3}(A^{\ast})^{2}+\frac{1}{192}iA^{3}+\frac{1}{16}%
iA^{2}A^{\ast}+\frac{1}{48}iA^{2}(A^{\ast})^{3}-\frac{1}{16}iA
\end{array}
\right)  \left(  t-t_{0}\right)  ^{2}e^{i\left(  t-t_{0}\right)  }\bigskip
$\newline$\left.  {}\right.  \hspace{0.5in}+\left(  -\frac{5}{16}A^{4}%
(A^{\ast})^{3}+\frac{9}{16}A^{3}(A^{\ast})^{2}-\frac{13}{48}A^{2}A^{\ast
}+\frac{1}{48}A\right)  \left(  t-t_{0}\right)  ^{3}e^{i\left(  t-t_{0}%
\right)  }\bigskip$\newline$\left.  {}\right.  \hspace{0.5in}+\left(
\begin{array}
[c]{c}%
+\frac{1}{4608}iA^{7}+\frac{7}{512}iA^{6}A^{\ast}-\frac{37}{1536}%
iA^{5}(A^{\ast})^{2}-\frac{1}{128}iA^{5}-\frac{1}{512}iA^{4}(A^{\ast}%
)^{3}\bigskip\\
+\frac{1}{128}iA^{3}+\frac{7}{256}iA^{4}A^{\ast}%
\end{array}
\right)  e^{3i\left(  t-t_{0}\right)  }\bigskip$\newline$\left.  {}\right.
\hspace{0.5in}+\left(
\begin{array}
[c]{c}%
-\frac{1}{96}A^{6}A^{\ast}-\frac{7}{64}A^{5}(A^{\ast})^{2}+\frac{1}{128}%
A^{5}+\frac{1}{384}A^{4}(A^{\ast})^{3}\bigskip\\
+\frac{21}{128}A^{4}A^{\ast}-\frac{1}{16}A^{3}%
\end{array}
\right)  \left(  t-t_{0}\right)  e^{3i\left(  t-t_{0}\right)  }\bigskip
$\newline$\left.  {}\right.  \hspace{0.5in}+\left(  -\frac{5}{64}%
iA^{5}(A^{\ast})^{2}+\frac{1}{8}iA^{4}A^{\ast}-\frac{3}{64}iA^{3}\right)
\left(  t-t_{0}\right)  ^{2}e^{3i\left(  t-t_{0}\right)  }\bigskip$%
\newline$\left.  {}\right.  \hspace{0.5in}+\left(  \frac{17}{2304}iA^{5}%
-\frac{17}{1536}iA^{6}A^{\ast}-\frac{5}{4608}iA^{7}\right)  e^{5i\left(
t-t_{0}\right)  }+\left(  \frac{5}{384}A^{6}A^{\ast}-\frac{5}{384}%
A^{5}\right)  \left(  t-t_{0}\right)  e^{5i\left(  t-t_{0}\right)  }\bigskip
$\newline$\left.  {}\right.  \hspace{0.5in}+\frac{1}{1152}iA^{7}e^{7i\left(
t-t_{0}\right)  }+c.c.$

We choose the homogeneous parts to the solutions $y_{1}$, $y_{2}$ and $y_{3}$
in such a manner that the solutions vanish at the initial time $t_{0}$, i.e.
$y_{1}\left(  t_{0}\right)  =y_{2}\left(  t_{0}\right)  =y_{3}\left(
t_{0}\right)  =0$. Next, we renormalize the integration constant $A$ and
create a new renormalized quantity $\mathcal{A}$ as%
\[
A=\mathcal{A}+a_{1}\varepsilon+a_{2}\varepsilon^{2}+a_{3}\varepsilon
^{3}+O\left(  \varepsilon^{4}\right)
\]
where the coefficients $a_{1},a_{2},a_{3},\ldots$ are chosen to absorb the
homogeneous parts of the solutions $y_{1},y_{2},\ldots$. Choosing%
\begin{align*}
a_{1}  &  =-\frac{i}{24}\mathcal{A}^{3},~a_{2}=-\frac{\mathcal{A}^{3}}%
{32}\left(  1-\frac{3}{2}\mathcal{AA}^{\ast}+\frac{1}{6}\mathcal{A}%
^{2}\right)  ,\\
a_{3}  &  =\frac{1}{1152}i\mathcal{A}^{7}-\frac{17}{1536}i\mathcal{A}%
^{6}\mathcal{A}^{\ast}+\frac{17}{2304}i\mathcal{A}^{5}-\frac{37}%
{1536}i\mathcal{A}^{5}(\mathcal{A}^{\ast})^{2}+\frac{7}{256}i\mathcal{A}%
^{4}\mathcal{A}^{\ast}+\frac{1}{128}i\mathcal{A}^{3}%
\end{align*}
we obtain\newline$y_{0}\left(  t\right)  =\mathcal{A}e^{i\left(
t-t_{0}\right)  }+c.c.\bigskip\newline y_{1}\left(  t\right)  =\left(
\frac{1}{2}\mathcal{A}\left(  1-\mathcal{AA}^{\ast}\right)  \left(
t-t_{0}\right)  e^{i\left(  t-t_{0}\right)  }-\frac{1}{24}i\mathcal{A}%
^{3}e^{3i\left(  t-t_{0}\right)  }\right)  +c.c.\bigskip\newline y_{2}\left(
t\right)  =\left(  \frac{1}{16}i\mathcal{A}^{3}(\mathcal{A}^{\ast})^{2}%
-\frac{1}{8}i\mathcal{A}\right)  \left(  t-t_{0}\right)  e^{i\left(
t-t_{0}\right)  }+\frac{1}{8}\mathcal{A}\left(  \mathcal{AA}^{\ast}-1\right)
\left(  3\mathcal{AA}^{\ast}-1\right)  ~\left(  t-t_{0}\right)  ^{2}%
~e^{i\left(  t-t_{0}\right)  }\bigskip\newline\left.  {}\right.
\hspace{0.5in}+\left(  \frac{3}{64}\mathcal{A}^{4}\mathcal{A}^{\ast}-\frac
{1}{32}\mathcal{A}^{3}\right)  e^{3i\left(  t-t_{0}\right)  }+\frac{1}%
{16}i\mathcal{A}^{3}\left(  \mathcal{AA}^{\ast}-1\right)  ~\left(
t-t_{0}\right)  ~e^{3i\left(  t-t_{0}\right)  }\bigskip\newline\left.
{}\right.  \hspace{0.5in}-\frac{1}{192}\mathcal{A}^{5}e^{5i\left(
t-t_{0}\right)  }+c.c.\bigskip$\newline$y_{3}\left(  t\right)  =\left(
-\frac{13}{128}\mathcal{A}^{4}(\mathcal{A}^{\ast})^{3}+\frac{11}%
{64}\mathcal{A}^{3}(\mathcal{A}^{\ast})^{2}\right)  \left(  t-t_{0}\right)
e^{i\left(  t-t_{0}\right)  }\bigskip\newline\left.  {}\right.  \hspace
{0.5in}+\left(  -\frac{3}{32}i\mathcal{A}^{4}(\mathcal{A}^{\ast})^{3}+\frac
{3}{32}i\mathcal{A}^{3}(\mathcal{A}^{\ast})^{2}+\frac{1}{16}i\mathcal{A}%
^{2}\mathcal{A}^{\ast}-\frac{1}{16}i\mathcal{A}\right)  \left(  t-t_{0}%
\right)  ^{2}e^{i\left(  t-t_{0}\right)  }\bigskip\newline\left.  {}\right.
\hspace{0.5in}+\left(  -\frac{5}{16}\mathcal{A}^{4}(\mathcal{A}^{\ast}%
)^{3}+\frac{9}{16}\mathcal{A}^{3}(\mathcal{A}^{\ast})^{2}-\frac{13}%
{48}\mathcal{A}^{2}\mathcal{A}^{\ast}+\frac{1}{48}\mathcal{A}\right)  \left(
t-t_{0}\right)  ^{3}e^{i\left(  t-t_{0}\right)  }\bigskip\newline\left.
{}\right.  \hspace{0.5in}+\left(  -\frac{37}{1536}i\mathcal{A}^{5}%
(\mathcal{A}^{\ast})^{2}+\frac{1}{128}i\mathcal{A}^{3}+\frac{7}{256}%
i\mathcal{A}^{4}\mathcal{A}^{\ast}\right)  e^{3i\left(  t-t_{0}\right)
}\bigskip\newline\left.  {}\right.  \hspace{0.5in}+\left(  -\frac{7}%
{64}\mathcal{A}^{5}(\mathcal{A}^{\ast})^{2}+\frac{21}{128}\mathcal{A}%
^{4}\mathcal{A}^{\ast}-\frac{1}{16}\mathcal{A}^{3}\right)  \left(
t-t_{0}\right)  e^{3i\left(  t-t_{0}\right)  }\bigskip\newline\left.
{}\right.  \hspace{0.5in}+\left(  -\frac{5}{64}i\mathcal{A}^{5}(\mathcal{A}%
^{\ast})^{2}+\frac{1}{8}i\mathcal{A}^{4}\mathcal{A}^{\ast}-\frac{3}%
{64}i\mathcal{A}^{3}\right)  \left(  t-t_{0}\right)  ^{2}e^{3i\left(
t-t_{0}\right)  }\bigskip\newline\left.  {}\right.  \hspace{0.5in}+\left(
\frac{17}{2304}i\mathcal{A}^{5}-\frac{17}{1536}i\mathcal{A}^{6}\mathcal{A}%
^{\ast}\right)  e^{5i\left(  t-t_{0}\right)  }+\left(  \frac{5}{384}%
\mathcal{A}^{6}\mathcal{A}^{\ast}-\frac{5}{384}\mathcal{A}^{5}\right)  \left(
t-t_{0}\right)  e^{5i\left(  t-t_{0}\right)  }\bigskip\newline\left.
{}\right.  \hspace{0.5in}+\frac{1}{1152}i\mathcal{A}^{7}e^{7i\left(
t-t_{0}\right)  }+c.c.$

\begin{remark}
DeVille \cite{DeVille RG} have obtained same result correct upto $O\left(
\varepsilon^{3}\right)  $. However, their computed expression of $a_{2}$ is
not correct. We have made the correction in the expression of $a_{2}$.
\end{remark}

We observe that each of $y_{1}\left(  t\right)  $, $y_{2}\left(  t\right)  $,
$y_{3}\left(  t\right)  $ contains secular terms. As a consequence the
solution%
\[
y\left(  t\right)  =y_{0}\left(  t\right)  +y_{1}\left(  t\right)
\varepsilon+y_{2}\left(  t\right)  \varepsilon^{2}+y_{3}\left(  t\right)
\varepsilon^{3}+O\left(  \varepsilon^{4}\right)
\]
becomes divergent as $t\rightarrow\infty$. To regularize the perturbation
series using RGM an arbitrary time $\tau$ is introduced and $t-t_{0}$ is split
as $\left(  t-\tau\right)  +\left(  \tau-t_{0}\right)  $. The terms containing
$\tau-t_{0}$ is observed in the renormalized counterpart $\mathcal{A}$ of the
constant of integration $A$. Since the final solution should not depend upon
the choice of the arbitrary time $\tau$, so%
\begin{equation}
\left.  \frac{\partial y}{\partial\tau}\right\vert _{\tau=t}=0
\label{CGO RG Cond}%
\end{equation}
for any $t$. However, DeVille et. al. \cite{DeVille RG} have simplified this
condition and proposed an equivalent condition as%
\begin{equation}
\left.  \frac{\partial y}{\partial t_{0}}\right\vert _{t_{0}=t}=0
\label{DeVille RG Cond}%
\end{equation}
We note that renormalized counterpart $\mathcal{A}$ is no longer a constant of
motion in RGM. The RG condition $\left(  \text{\ref{DeVille RG Cond}}\right)
$ is developed in such a manner that one need to differentiate the terms
containing $e^{i\left(  t-t_{0}\right)  }$, $e^{-i\left(  t-t_{0}\right)  }$,
$\left(  t-t_{0}\right)  e^{i\left(  t-t_{0}\right)  }$ and $\left(
t-t_{0}\right)  e^{-i\left(  t-t_{0}\right)  }$ and thereafter substituting
$t_{0}=t$ the resultant expression is equated to zero. The other terms related
to higher harmonics are not involved in RG condition. Simplifying RG condition
$\left(  \text{\ref{DeVille RG Cond}}\right)  $ we get%
\[
\frac{\partial\mathcal{A}}{\partial t_{0}}=\mathcal{A}i-\frac{1}{2}%
\mathcal{A}\left(  \mathcal{AA}^{\ast}-1\right)  \varepsilon-\frac{1}%
{8}i\mathcal{A}\left(  1-\frac{1}{2}\mathcal{A}^{2}(\mathcal{A}^{\ast}%
)^{2}\right)  \varepsilon^{2}-\frac{1}{64}\mathcal{A}^{3}(\mathcal{A}^{\ast
})^{2}\left(  \frac{13}{2}\mathcal{AA}^{\ast}-11\right)  \varepsilon^{3}
\]
to the order $O\left(  \varepsilon^{4}\right)  $. Taking $\mathcal{A=}\frac
{R}{2}e^{i(t+\theta)}$ we obtain corresponding amplitude and phase flow
equations to the order $O\left(  \varepsilon^{4}\right)  $ as%
\begin{align}
\frac{dR}{dt}  &  =\frac{1}{2}R\left(  1-\frac{R^{2}}{4}\right)
\varepsilon+\frac{1}{1024}R^{5}\left(  11-\frac{13}{8}R^{2}\right)
\varepsilon^{3}+O\left(  \varepsilon^{4}\right) \label{Amp Eq Order 3}\\
\frac{d\theta}{dt}  &  =-\frac{1}{8}\left(  1-\frac{R^{4}}{32}\right)
\varepsilon^{2}+O\left(  \varepsilon^{4}\right)  \label{Phase Eq Order 3}%
\end{align}
To the authors' knowledge these higher order flow equations are reported for
the first time in the literature. We remark that above flow equations match
exactly with $O(\varepsilon^{3})$ flow equations of the Van der Pol equation
\cite{Sarkar Bhattacharjee}. Although not done explicitly, we expect that the
$O\left(  \varepsilon^{4}\right)  $ VdP flow equations would also have the
equivalent forms. For latter reference, we also write down the order
$O(\varepsilon^{2})$ solution of the Rayleigh equation \cite{Chen G O RG}
\begin{equation}
\label{Sol1}y(t)= R(t)\cos(t+\theta) + \frac{\varepsilon}{96}R(t)^{3}%
(\sin3(t+\theta)-\sin(t+\theta))
\end{equation}

\begin{figure}[h]
\begin{center}%
\begin{tabular}
[c]{cc}%
\includegraphics[width=3in]{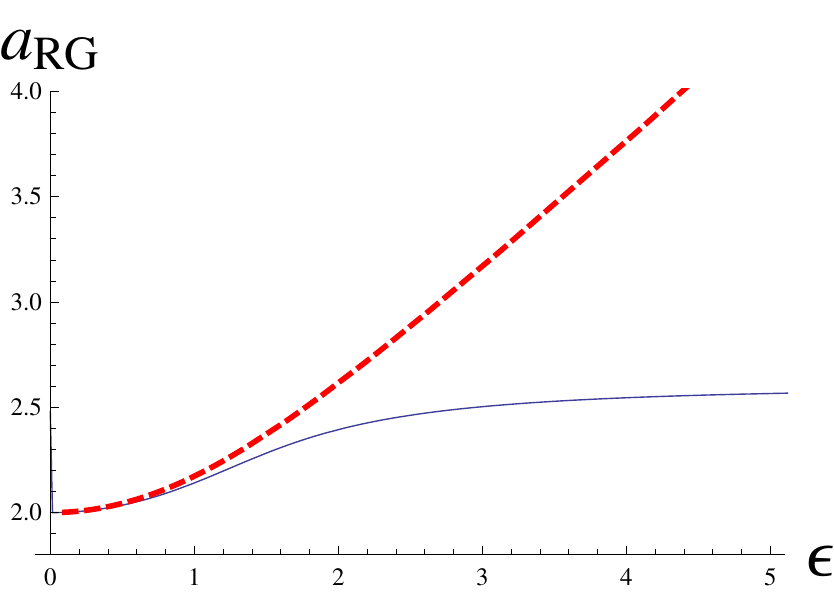} &
\includegraphics[
width=3in]{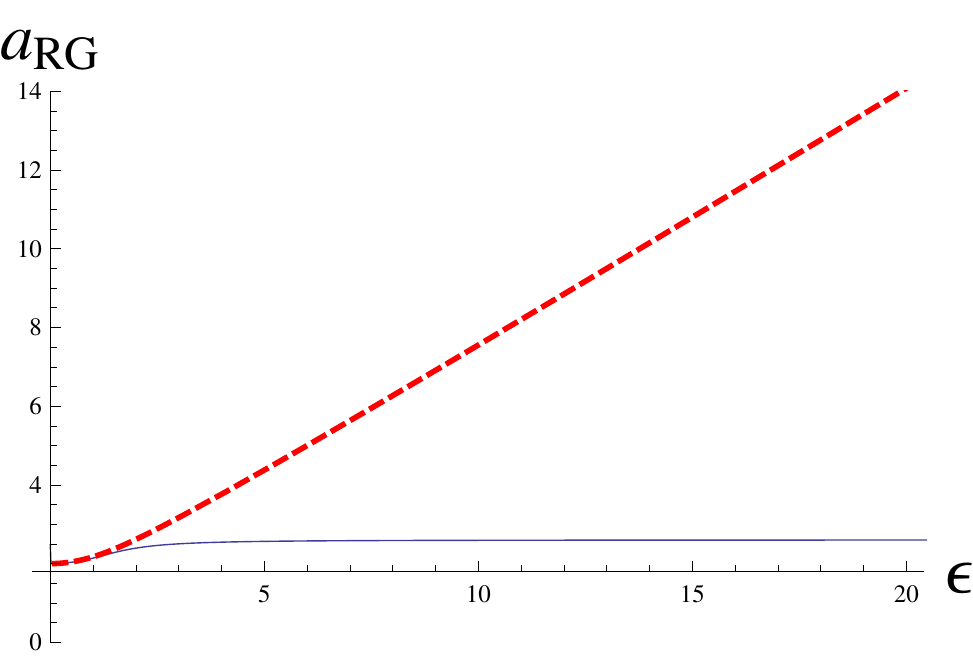}\\
$\left(  a\right)  $ & $\left(  b\right)  $%
\end{tabular}
\end{center}
\caption{Graph of $a_{RG}\left(  \varepsilon\right)  $ $($by solid lines$)$
correct upto $O\left(  \varepsilon^{4}\right)  $ and compared with exact graph
of $a\left(  \varepsilon\right)  $ $($by dotted lines$)$ for $0<\varepsilon
\leq5$ in $\left(  a\right)  $ and for $0<\varepsilon\leq20$ in $\left(
b\right)  $.}%
\label{Fig RG O(3) Compare}%
\end{figure}

Solving the amplitude equation $\left(  \text{\ref{Amp Eq Order 3}}\right)  $
by numerical method and taking the limit as $t\rightarrow\infty$ so that for a
fixed value of $\varepsilon$ we have $R\rightarrow a_{RG}\left(
\varepsilon\right)  $, the approximation of the amplitude of limit cycle of
Rayleigh equation $\left(  \text{\ref{Rayleigh}}\right)  $ by RGM, we obtain
Figure \ref{Fig RG O(3) Compare} representing $\varepsilon$ dependence of the
amplitude $a_{RG}$ by solid lines. 

Thus we observe that the RG flow equation to the order $O\left(
\varepsilon^{4}\right)  $ for the amplitude does not give good approximation
to the exact solution for moderate and large values of $\varepsilon$.

\section{Improved RG Method: Nonlinear Time\label{Sec M-RG Sol Description}}

In RGM an arbitrary time $\tau$ is introduced in between current time $t$ and
the initial time $t_{0}$ so that $t-t_{0}=\left(  t-\tau\right)  +\left(
\tau-t_{0}\right)  $ in order to remove the divergent terms in the naive
perturbation expansion for the solution of the given differential equation.
The solution is renormalized by suitable choice of the constants of
integration to remove the terms containing $\left(  \tau-t_{0}\right)  $ and
keeping the terms having $\left(  t-\tau\right)  $. Since the solution should
be independent of the arbitrary time $\tau$, the RG condition%
\[
\left.  \frac{\partial y}{\partial\tau}\right\vert _{\tau=t}=0
\]
is applied to the renormalized solution. However, in the previous section we
have seen that the method fails to produce good approximations to the exact
solution for $\varepsilon\sim O(1)$. Our target is not only to remove the
divergent terms in the solution but also to introduce some control parameter
$h(\varepsilon)$ which can control the RG solution in such a manner that this
solution ultimately converges to the exact solution. Moreover, our another
goal is to achieve this accuracy by merely solving the differential equation
to a minimal order of the expansion parameter, viz., upto $O\left(
\varepsilon^{2}\right)  $ or less.

Since the basic idea is to split the time difference $t-t_{0}$ by introduction
of an arbitrary time, so we can write $t-t_{0}=\left(  t-\dfrac{\tau
}{\varepsilon}\right)  +\left(  \dfrac{\tau}{\varepsilon}-t_{0}\right)  $.
From now on let us assume that $0<<\varepsilon<\approx1$. The case
$\varepsilon>\approx1$ will be commented upon later. The constants of
integration can be renormalized in order to remove the terms containing
$\left(  \dfrac{\tau}{\varepsilon}-t_{0}\right)  $ from the solution keeping
the terms containing $\left(  t-\dfrac{\tau}{\varepsilon}\right)  $. Finally
analogous to the classical RG method we put $t=\dfrac{\tau}{\varepsilon}$,
i.e. $\tau=\varepsilon t$, in%
\begin{equation}
\frac{\partial y}{\partial\tau}=0 \label{Pre RG eq}%
\end{equation}
giving rise to an improved form of the RG flow equation to remove secular
terms involving $\left(  t-\dfrac{\tau}{\varepsilon}\right)  $. So far the
improved method does not produce any qualitative new result compared to the
RGM and so we must get the same phase and amplitude equation as deduced in
Section \ref{Sec RG Sol}.

We next proceed one step further. As stated already in the Introduction, we
now exploit the possibility of extending the original linear $t$ dependence of
$\tau$ viz., $\tau=t$ of RGM in removing the explicit divergences by a
\emph{nonlinear dependence} $\tau=\varepsilon t$ along with the
\emph{additional condition} that $\tau\rightarrow\varepsilon^{-n}\phi
(\tilde\tau)$, where $\phi$ a slowly varying scaling function of the O(1)
rescaled variable $\tilde\tau=\varepsilon\tilde t\sim O(1)$, as the original
linear time $t\rightarrow\infty$ following the\emph{ scales} $t\sim
\varepsilon^{-n}\tilde t,\ n=1,2,\ldots$. (Note that linear time flows with
uniform rate 1 and $\tau$ is nonlinear since rate $\dot\phi(\tilde\tau)<1$).
It follows that for a given nonlinear differential system, such a nonlinear
time dependence always exists and nontrivial, provided one invokes a
\emph{duality principle} transferring nonlinear influences from the far
asymptotic region into the finite observable sector in a {cooperative }manner
\cite{DPD, DPD Sen, DPD New}. In Appendix, we give a brief overview of the
\emph{novel} analytic framework extending the standard classical analysis to
one that supports naturally the above stated \emph{duality structure} and the
emergent nonlinear scaling patterns typical for a given nonlinear system.

In fact, as the linear time $t\rightarrow\infty$ following the above hierarchy
of scales, there exists $\tilde{t}_{n}$ such that $1<<\left(  \varepsilon
t\right)  ^{n}<\varepsilon^{-n}<\tilde{t}_{n}$ and satisfying \emph{the
inversion law} $\tilde{t}_{n}/\varepsilon^{-n}\propto\varepsilon^{-n}/\left(
\varepsilon t\right)  ^{n}$. This inversion law makes a room for transfer of
effective influences, typical for the nonlinear system concerned, from
nonobservable sector $t>\varepsilon^{-n}$ to the observable sector
$t<\varepsilon^{-n}$ bypassing the dynamically generated singular points
denoted by the scales $\varepsilon^{-n}$. Notice the nonlinear connection
between scales of the form $\varepsilon^{n} \tilde t_{n}$ with the scale
$\varepsilon t$ via duality structure (c.f. Appendix). Let $\tilde{t}\left(
t\right)  =\underset{n\rightarrow\infty}{\lim}\left(  \tilde{t}_{n}\right)
^{1/n}$ so that $\varepsilon t<\varepsilon^{-1}<\tilde{t}\left(  t\right)  $
and $\tilde{t}/\varepsilon^{-1}\propto\varepsilon^{-1}/\left(  \varepsilon
t\right)  $. Define
\begin{equation}
\label{visible}h_{0}\left(  \tilde\tau\right)  =\underset{n\rightarrow\infty
}{\lim}\log_{\varepsilon^{-n}}\tilde{t}_{n}/\varepsilon^{-n}.
\end{equation}
Here, the scaling exponent $h_{0}$ corresponds to the visibility norm
(Appendix), that can access (encode) the non-perturbative region (information)
of the nonlinear system and $\tilde\tau$ is an O(1) rescaled variable. The
exponent $h_{0}$ is \emph{scaling invariant} in the sense that it appears
uniformly for every $n$ as $t\rightarrow\infty$ through the scales $\tilde
t_{n}=\varepsilon^{-n} \varepsilon^{-n h_{0}(\tilde\tau)}$. As a consequence,
a significant amount of asymptotic scaling information in the limit
$t\rightarrow\infty$ could be simply retrieved by considering the scaling
limit instead at $t=\varepsilon^{-1}$.

Exploiting the above insight, one now writes the nonperturbative scaling limit
in the form
\begin{equation}
\tau=\lim\varepsilon t=\varepsilon^{-h_{RG}\left(  \tilde\tau\right)
}>1,\ \varepsilon<1 \label{scaling1}%
\end{equation}
as $t\rightarrow\varepsilon^{-1}$. Moreover, $h_{RG}\left(  \tilde\tau\right)
=1-h_{0}\left(  \tilde\tau\right)  $. As noted already, the scaling exponent
$h_{0}\left(  \tilde\tau\right)  $ here encodes the \emph{effective
}cooperative influence of far asymptotic sector $t>\varepsilon^{-n}$ into the
observable sector $1<t<\varepsilon^{-n}$ by the inversion mediated duality
principle. As pointed out in Appendix, the duality principle \emph{does} allow
asymptotic limiting (non-perturbative) behaviour of the nonlinear system to be
encoded into the scaling exponents of the nonlinear time $\tau$ that, in turn,
offers an efficient handle in uncovering key dynamical information of the said
system. Notice that, in the absence of the said duality the linear time $t$
can in principle attain the scale $\varepsilon^{-1}$ $\left(  \text{say}%
\right)  $, and as a consequence $h_{0}=0$, retrieving the ordinary scaling of
$\tau=\varepsilon t\sim\varepsilon^{-1}$ as $t\sim\varepsilon^{-2}$. This also
establishes, in retrospect, that the scaling exponent $h_{0}\left(
\varepsilon\right)  $ is well defined and can exist nontrivially i.e.
$h_{0}\sim O\left(  1\right)  $ in a nonlinear problem (c.f. Appendix). As a
consequence, \emph{the RG control parameter }$\emph{h}_{RG}$\emph{ can be of
both the signs, with relatively small numerical value in fully developed
nonlinear systems }$\varepsilon\gg1$\emph{, but with a possible }$O\left(
1\right)  $\emph{ variations for }$\varepsilon\sim O\left(  1\right)  $\emph{
or less.}

The above construction actually tells somewhat more. Corresponding to the
first generation scales $\varepsilon^{-n}$, one can, in fact, have the
\emph{second generation } nonlinear scales
\begin{equation}
\tau_{m}=\lim\varepsilon^{m} t=\varepsilon^{-m h^{m}_{RG}\left(  \tilde
\tau\right)  }>1,\ \varepsilon<1, \ m>1 \label{scaling2}%
\end{equation}
as $t\rightarrow\varepsilon^{-m}$ with $h^{1}_{RG}=h_{RG}$. The
\emph{nonlinear time} $\tau$ now stands for these hierarchy of scales
$\{\tau_{m}\}$. Consequently, as the linear time $t$ approaches $\infty$
through the first generation linear scales, the slowly varying nonlinear time
$\tau$ approaches either to $\infty$ or 0 at slower and slower rates as
represented by the numerically small RG scaling exponents $h^{m}%
_{RG}(\varepsilon)$, each of which remains almost constant over longer and
longer intervals of $\varepsilon^{-1}$ (as $\varepsilon^{-1} \rightarrow
\infty$ ). In the present paper we show how the first two scaling exponents
$h_{i}(\tilde\tau), \ i=1,2$ relate to the nonperturbative properties of the
limit cycle. We expect higher order scaling exponents $h_{m}$ would have vital
role in bifurcation of nonautonomous systems. This problem will be
investigated elsewhere.

Let us remark that for $\varepsilon>1$, we consider instead the first
generation scales as $\varepsilon^{n}$, and the duality is invoked for
variables satisfying $t/\varepsilon<\varepsilon<\tilde{t}(t)$ so that the
asymptotic scaling variables are derived as $\tau_{m}=\varepsilon^{m
h^{m}_{RG}(\tilde\tau)},\ \varepsilon>1$ where $h^{m}_{RG}=1-h^{m}_{0}$.
Moreover, said proliferation of nonlinear scales (\ref{scaling2}) actually
continues ad infinitum. In fact, interpreting each second generation scale
$\tau_{m}$, $m$ fixed, as first generation scale, and iterating above steps
one associates third generation scales $\tau_{m_{k}}, \ k=1,2,\ldots$, and so on.

It now follows, from the above general remarks on the behaviour of $h_{RG}$,
that the nonlinear time $\tau$ actually approaches 0 or $\infty$ as $\tau
\sim(\log\varepsilon)^{-\alpha}$ or $\tau\sim(\log\varepsilon)^{\alpha},
\ \alpha>0$ respectively as $\varepsilon\rightarrow\infty$. However, one must
have $\tau=\varepsilon^{-h_{RG}(\varepsilon)}\rightarrow\infty$ as
$\varepsilon\rightarrow0$. An example of the asymptotic behaviours of $h_{RG}$
is given by $\tau_{m}=\varepsilon^{\pm\alpha_{m}\frac{\log\log\varepsilon
}{\log\varepsilon}}$ for $\varepsilon\rightarrow\infty$, which one expects to
verify explicitly in evaluation of asymptotic quantities, such as amplitude of
a periodic cycle, in a nonlinear system.

In the IRGM, we exploit this duality induced nontrivial scaling information to
rewrite the lowest order perturbative flow equations (\ref{Amp Eq Order 3})
and (\ref{Phase Eq Order 3}) as the asymptotic RG flow equations in the limit
$t\rightarrow\infty$
\begin{align}
\frac{da}{d\tau_{1}} &  =\frac{1}{2}a\left(  1-\frac{a^{2}}{4}\right)
\label{Amp}\\
\frac{d\psi}{d\tau_{2}} &  =-\frac{1}{8}\left(  1-\frac{a^{4}}{32}\right)
\label{Phase}%
\end{align}
for the amplitude $a=a(\tilde{\tau})$ and the phase $\psi=\psi(\tilde{\tau})$
of the limit cycle of both the Rayleigh and Van der Pol equations, involving
slowly varying nonlinear time scales $\tau_{i},\ i=1,2$. The asymptotic
scaling functions $\tau_{1}=\phi_{1}(\varepsilon t)=\varepsilon^{h_{RG}^{1}}$
and $\tau_{2}=\phi(\varepsilon^{2}t)=\varepsilon^{2h_{RG}^{2}}$ are activated
invoking nonlinear limits as in (\ref{scaling2}) as $t\rightarrow\varepsilon
^{-1}$ and $t\rightarrow\varepsilon^{-2}$ successively in the above equations.
The slowly varying almost constant scaling functions $\phi_{1}$ and $\phi_{2}%
$, satisfying $|\ddot{\phi}_{i}|<<|\dot{\phi}_{i}^{2}|<<1$, are assumed to
have a rhythmic pattern over the cycle: when $\phi_{1}$ varies slowly,
$\phi_{2}$ remains almost constant i.e. $\dot{\phi}_{1}>0,\ \dot{\phi}%
_{2}\approx0$ and vice versa successively on the cycle (c.f. Appendix Sec. B).
Nontrivial ultrametric neighbourhood structure induced asymptotically by
duality principle (c.f. Appendix Sec. A) can indeed support such \emph{locally
constant} nonlinear rhythmic behaviour. The above flow equations may therefore
be considered \emph{exact} and encode  non-perturbative information of
the limit cycle variables $a$ and $\psi$ respectively. The conventional
perturbative RG flow equations in the linear time $t$ is now extended into the
non-perturbative flow equations in the nontrivial scaling variable $\tau
_{i}=\varepsilon^{ih_{RG}^{i}(\tilde{\tau})},\ i=1,2$ involving the
nonlinearity parameter $\varepsilon>1$. The perturbative fixed point for the
amplitude equation at $a=2$ for $t\rightarrow\infty$ corresponding to the
periodic oscillation with $\varepsilon<<1$ is superimposed by \emph{small scale
periodic} flow of amplitude $a(\tau_{1})$ over the entire cycle. The
associated phase $\psi(\tau_{2})$ then flow at a slower rate \emph{linearly
with the higher order scale $\tau_{2}$} when $a(\tau_{1})$ remains almost
constant over a relatively small period of time.

The RG estimated approximate formulae for the amplitude $a(\varepsilon)$ for
the Rayleigh and Van der Pol limit cycles are obtained from the equation
(\ref{Amp}) in the Sec.4.1, when appropriate boundary condition, derived
either from exact computation or from perturbative analysis, is used for a
suitable finite value of $\varepsilon$. In the next subsection 4.2, we present
the efficient graphs of the Rayleigh and VdP limit cycle parametrized by the
nonlinear scales $\tau_{i}=\phi_{i}(\tilde\tau), \ \tilde\tau\sim O(1)$ for
fixed values of the nonlinearity parameter $\varepsilon$.

As it turns out, entire onus in the improved RG analysis essentially rests in
proper estimation/identification of the scaling functions $h^{i}_{RG}$ (i.e.
$\phi_{i}(\tilde\tau)$) which should yield correct dynamical properties of a
nonlinear system. We hope to undertake more detailed and systematic analysis
for determining $h^{i}_{RG}$ elsewhere. In this work we limit ourselves only
to show that IRGM can indeed yield correct amplitude and solution for the
Rayleigh and VdP systems provided one makes appropriate choice of $h^i_{RG}$ based on
clues from exact computations and previously known approximate results (for
instance the perturbative RGM). We remark finally that the perturbative RG method is known 
to extend the conventional multiple scale method \cite{Chen G O RG}. Nonlinear time formalism 
introduces new set of nonlinear scales $\phi_n(\tilde \tau)$ associated 
with ordinary scales $\varepsilon^n$. 
We study here the nontrivial applications of such nonlinear scaling functions.

\subsection{Approximate Formula for
Amplitude\label{SubSec Amplitude Estimation}}

We shall now use the above asymptotic amplitude flow equation (\ref{Amp}) to
find analytic approximations of the amplitudes of the limit cycle for both the
Rayleigh and Van der Pol equations.

\begin{figure}[h]
\begin{center}%
\begin{tabular}
[c]{cc}%
\includegraphics[width=3in]{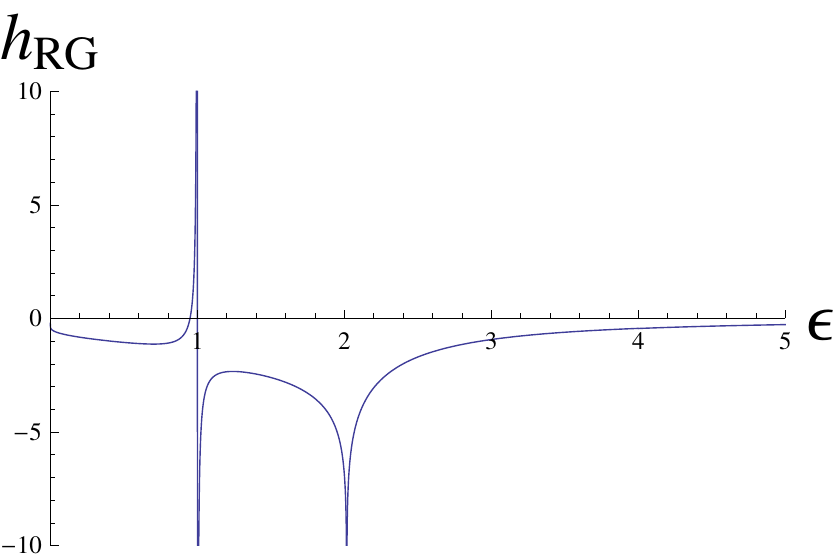} &
\includegraphics[
width=3in]{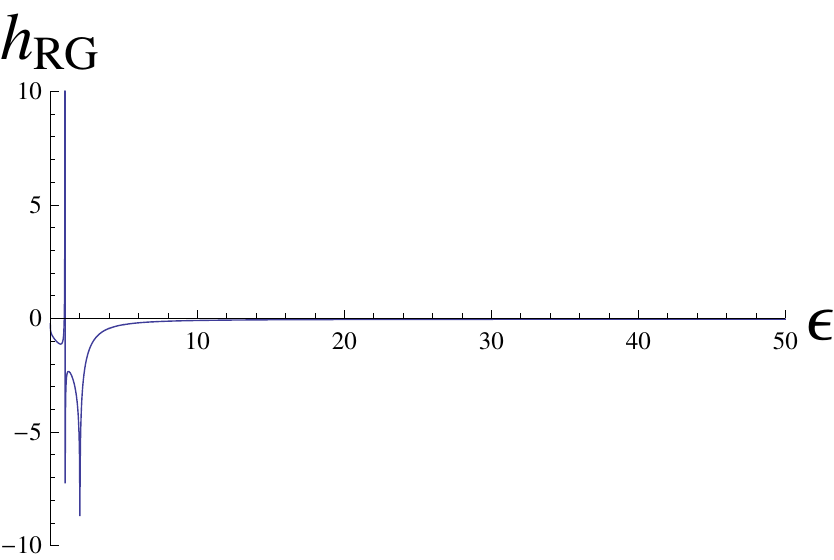}\\
$\left(  a\right)  $ & $\left(  b\right)  $%
\end{tabular}
\end{center}
\caption{The graph of $h_{RG}\left(  \varepsilon\right)  $ used for
approximation of the amplitude of the Rayleigh equation $\left(
\text{\ref{Rayleigh}}\right)  $ by HAM given by $\left(
\text{\ref{Cont Parameter RG}}\right)  $ for $0<\varepsilon\leq5$ in $\left(
a\right)  $ and for $0<\varepsilon\leq50$ in $\left(  b\right)  $.}%
\label{Fig h IRG Rayleigh(epsilon)}%
\end{figure}

By a direct integration, one obtains from (\ref{Amp})%

\begin{equation}
\ln\left(  a^{2}-4\right)  -2\ln a=-\varepsilon^{h_{RG}}%
-0.87953\label{Amp Sol RG RL}%
\end{equation}
as the Rayleigh limit cycle amplitude where we use the boundary condition the
value $a=2.17271$ for $\varepsilon=1$ (this choice simplifies calculation). It
follows immediately that for suitable choices of the control parameter
$h_{RG}$ one can achieve efficient matching for the estimated amplitude
$a_{{E}}(\varepsilon)$. For example, using the HAM generated approximate
formula (\ref{HAM Amp New}) for $a_{{E}}(\varepsilon)$, we can determine the
control parameter $h_{RG}(\varepsilon)$ by the formula
\begin{equation}
h_{RG}=\frac{1}{\ln\varepsilon}\ln\left\{  \left\vert \ln\left(  \frac{a^{2}%
}{a^{2}-4}\right)  -0.87953\right\vert \right\}  \label{Cont Parameter RG}%
\end{equation}
In Figure \ref{Fig h IRG Rayleigh(epsilon)}, we display the typical piece-wise smooth form of $h_{RG}%
(\varepsilon)$ given by $\left(  \text{\ref{Cont Parameter RG}}\right)  $ for
the Rayleigh limit cycle amplitude that would reproduce the HAM generated
amplitude with relative error less that $1\%$. Clearly, the graph reveals
variability of $h_{RG}$ for moderate values of $\varepsilon$, but the
variability dies out fast for larger values $\varepsilon$, as expected.

\begin{figure}[h]
\begin{center}
\includegraphics[width=3in]{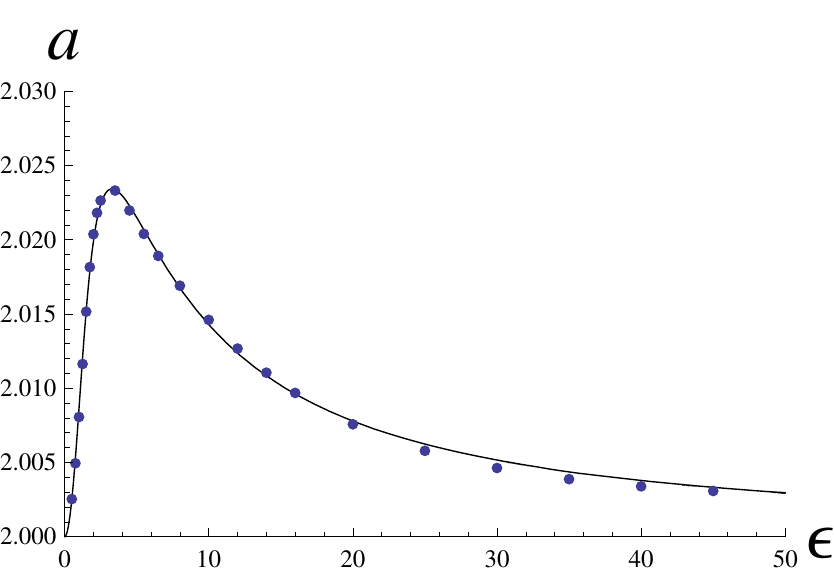}
\end{center}
\caption{The exact amplitude of Van der Pol Equation $\left(  \text{\ref{VdP}%
}\right)  $ $($by solid line$)$ and its approximation $a_{E}\left(
\varepsilon\right)  $ given by $\left(  \text{\ref{VdP Amplitude IRGM}%
}\right)  $ $($by bold points$)$ for $0<\varepsilon\leq50$.}%
\label{Fig VdP Amp}%
\end{figure}

We recall that the corresponding graph of the exact computed values of VdP
amplitude $a(\varepsilon)$, on the other hand, has a hump like shape with a
maximum roughly at $\varepsilon\approx2.0235$ and having the asymptotic limits
$2$ as $\varepsilon\rightarrow0$ and $\infty$. Lopez et al \cite{Lopez VDP
Amp} obtained HAM generated approximate formula for the VdP amplitude with
relative error less than $0.05\%$ at the order $O\left(  \varepsilon
^{4}\right)  $. It is interesting to note that the RG generated formula
$\left(  \text{\ref{Amp Sol RG RL}}\right)  $ can reproduce the exact computed
values of the VdP amplitude with error less than $0.05\%$ directly from only
the first order RG flow equation. To achieve this goal we first intuitively
guess a piecewise smooth formula for the estimated amplitude
$a_{E}$ by
\begin{equation}
a_{E}\left(  \varepsilon\right)  =\left\{
\begin{array}
[c]{lc}%
1.998+\dfrac{0.015}{8.121~e^{-2.139~\varepsilon}+0.512~e^{0.043~\varepsilon}}
& 0<\varepsilon<3\\
2.0025+\dfrac{0.031}{0.5~e^{-2.033\left(  \varepsilon-2.183\right)
}+1.869~e^{0.087\left(  \varepsilon-6.376\right)  }} & 3\leq\varepsilon\leq50
\end{array}
\right.  \label{VdP Amplitude IRGM}%
\end{equation}
keeping the maximum relative percentage error $\left\vert \dfrac{a_{E}\left(
\varepsilon\right)  -a\left(  \varepsilon\right)  }{a\left(  \varepsilon
\right)  }\times100\right\vert $ less than $0.05\%$. This shows that the
approximation is quite accurate. The graph of $a_{E}\left(  \varepsilon
\right)  $ is compared with the exact values in Figure \ref{Fig VdP Amp}. One
may as well use a least square fit of the exact data instead of the above fit.
We do not pursue this approach here. 

\begin{figure}[h]
\begin{center}%
\begin{tabular}
[c]{cc}%
\includegraphics[width=3in]{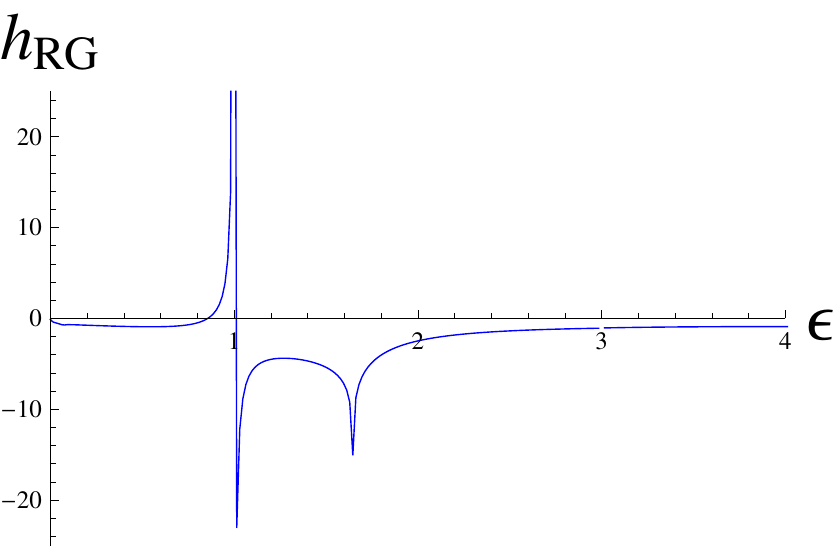} &
\includegraphics[
width=3in]{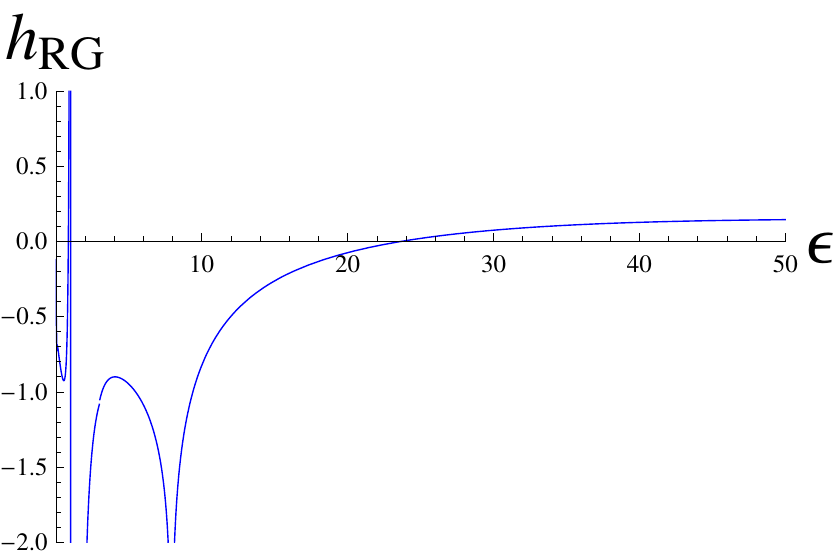}\\
$\left(  a\right)  $ & $\left(  b\right)  $%
\end{tabular}
\end{center}
\caption{The graph of $h_{RG}\left(  \varepsilon\right)  $ used for
approximation $\left(  \text{\ref{Cont Parameter RG2}}\right)  $ of the
amplitude of the Van der Pol equation $\left(  \text{\ref{VdP}}\right)  $ for
$0<\varepsilon\leq4$ in $\left(  a\right)  $ and for $0<\varepsilon\leq50$ in
$\left(  b\right)  $.}%
\label{Fig h IRG VDP(epsilon)}%
\end{figure}

Using this efficient formula for the VdP amplitude, we then obtain the RG flow
equation in the form
\begin{equation}
\ln\left(  a^{2}-4\right)  -2\ln a=-\varepsilon^{h_{RG}}-4.08785
\label{Amp Sol RG VDP}%
\end{equation}
where we use the boundary condition $a=2.0086$ for $\varepsilon=1$ (for
simplicity of calculation) for the VdP amplitude. Inverting this equation, we
finally obtain the corresponding RG control parameter
\begin{equation}
h_{RG}=\frac{1}{\ln\varepsilon}\ln\left\{  \left\vert \ln\left(  \frac
{a_{E}^{2}}{a_{E}^{2}-4}\right)  -4.08785\right\vert \right\}
\label{Cont Parameter RG2}%
\end{equation}
Figure \ref{Fig h IRG VDP(epsilon)} displays the piecewise smooth variation of $h_{RG}$ with
$\varepsilon$. The rapid $O\left(  1\right)  $ variation for moderate values
of $\varepsilon$ is evident in Fig.6(a). As expected, $h_{RG}$ dies out fast
for larger values of $\varepsilon$. However, a change in sign is noticed here
already for $\varepsilon>20$ $\left(  \text{Figure
\ref{Fig h IRG VDP(epsilon)}}\left(  \text{b}\right)  \right)  $. One expects
many more such small scale sign variations as $\varepsilon\rightarrow\infty$.
This particular form of the control parameter $h_{RG}$, in turn, would
reproduce the VdP amplitude with relative error less than $0.05\%$. As this
level of accuracy is achieved only at the order $O\left(  \varepsilon\right)
$, the improved RGM may be considered to be more efficient and advantageous
compared to the HAM.

Alternatively, the amplitude equation $\left(  \text{\ref{Amp}}\right)  $ can
be inverted as
\begin{equation}
a(\tilde\tau)=\frac{a_{0}}{\sqrt{e^{-\tilde\tau}+\frac{a_{0}^{2}}%
{4}(1-e^{-\tilde\tau})}} \label{amp}%
\end{equation}
where $a_{0}$ is estimated from the exact value of amplitude $a(\varepsilon
_{0})$ for a suitably chosen value of $\varepsilon$, for instance $\varepsilon=1$. 
Recall that for the VdP
equation $\tau\sim(\log\varepsilon)^{\alpha}$ and for the Rayleigh equation
$\tau\sim(\log\varepsilon)^{-\alpha}$ for $\varepsilon\rightarrow\infty$ and
$\alpha>0$. By adjusting suitably the values of $\alpha$ over appropriate
intervals on $\varepsilon$ one should be able to obtain efficient matchings
with the exact values of $a(\varepsilon)$.

To summarize, the recipe for deriving approximate formula for limit cycle
amplitude of a nonlinear system can be stated as follows: Determine the first
order (perturbative) RG flow equation for amplitude in the nonlinear time
$\tau$. This will yield an explicit formula for amplitude $a$ as a function of
the nonlinearity parameter $\varepsilon$ and the control parameter $h_{RG}$.
Efficient match with the exact amplitude can be achieved by right choice of
the control parameter $h_{RG}$ or $\alpha$. Alternatively, determine an
efficient formula for $a(\varepsilon)$ by inspection (expert guess) or by
appropriate curve fitting methods. Then determine the control parameter
$h_{RG}$ by an inversion of the estimated amplitude $a_{E}(\varepsilon)$ as in
equation $\left(  \text{\ref{Cont Parameter RG2}}\right)  $ (and Fig.6). Since
the equations concerned form a closed system, this already gives a (numerical)
proof of the unique existence of $h_{RG}$ for a given nonlinear oscillation.

\subsection{Approximate Limit Cycle}

Here we calculate the approximate limit cycle orbit for the Rayleigh and VdP
equations for a sufficiently large $\varepsilon>1$. Perturbative RGM fails to give 
correct relaxation oscillation solution. The first order solution given in Sec. 2 by 
HAM is also found insufficient. In this work we do not under take the problem of 
computing approximate limit cycle by HAM, which has been addressed by Lopez et al 
\cite{Lopez VDP Amp} for the VdP equation. Our aim here is to highlight the strength 
of IRGM over perturbative RGM.

For a sufficiently large
time $t\rightarrow\varepsilon^{n},\ n$ large, but fixed, the slowly varying
nonlinear scales $\tau_{1}$ and $\tau_{2}$ are activated in a successive
rhythmic manner, as explained in Sec. 4 (see also Appendix Sec.B), so that the
perturbative solution given in (\ref{Sol1}) is extended to the asymptotic
limit cycle (relaxation oscillation) solution
\begin{equation}
y(\tau_{1},\tau_{2})=a(\tau_{1})\cos(\varepsilon^{n}+\psi(\tau_{2}%
))+Y\label{Sol2}%
\end{equation}
where the amplitude $a$ and phase $\psi$ \emph{flow} along the cycle following
the nonperturbative flow equations (\ref{Amp}) and (\ref{Phase}) in the
asymptotic scaling variables $\tau_{1}$ and $\tau_{2}$ respectively. Here, $Y$
encodes all the renormalized perturbative terms depending on higher order,
slowly varying nonlinear scales $\tau_{i},\ i>2$. The corresponding velocity
component $\dot{y}=\frac{\partial y}{\partial t}$ at $t=\varepsilon^{n}$ has
the form
\begin{equation}
\dot{y}(\tau_{1},\tau_{2})=-a(\tau_{1})\sin(\varepsilon^{n}+\psi(\tau
_{2}))+\sum_{i}\dot{\tau}_{i}\frac{\partial \tilde Y}{\partial\tau_{i}}\label{Sol3}%
\end{equation}
where $\tilde Y$ represents possible slow variations in all the nonlinear scales 
$\tau_i, \ i\geq 1$.
Equations (\ref{Sol2}) and (\ref{Sol3}) are the parametric equations of the
limit cycle, parametrized by multiple nonlinear scales, when slowly varying
amplitude $a$ and phase $\psi$ are computed from (\ref{Amp}) and (\ref{Phase})
respectively. An alternative derivation of (\ref{Sol2}) and (\ref{Sol3}) based
purely on duality induced nonlinear scales is given in Appendix Sec. B. We
remark that (\ref{Sol2}) and (\ref{Sol3}) actually represent the general form
of the limit cycle for a much larger class of Lienard system having unique
limit cycle. \emph{Typical geometric shape of the periodic cycle of a given
nonlinear system is controlled entirely by the rhythmic cooperative, almost
constant variations of the nonlinear scales $\tau_{i}$.} As explained in Sec.
4 (See also Appendix Sec. A), scale invariance of the scaling functions
$\tau_{1}=\phi_{1}(\tilde{t}_{1}/\varepsilon)$ and $\tau_{2}=\phi_{1}%
(\tilde{t}_{2}/\varepsilon^{2})$ tells that as the linear time $t\rightarrow
\varepsilon^{n}$, $\tau_{1}\rightarrow\phi_{1}(1)=1$ and $\tau_{2}%
\rightarrow\phi_{2}(1)=1$ with $\tilde{t}_{1}\rightarrow\varepsilon$ and
$\tilde{t}_{2}\rightarrow\varepsilon^{2}$. We now set the initial conditions
$a(1)=a_{\mathrm{amp}}(\varepsilon),\ a^{\prime}(1)=\frac{a(1)}{2}%
(1-\frac{a(1)^{2}}{4})$ and $\psi(1)=0$, where $a_{\mathrm{amp}}(\varepsilon)$
is the exact (experimental) value of the amplitude of the limit cycle. Setting
further $\tau_{1}=1+\eta_{1}$ and $\tau_{2}=1+\eta_{2}$ as $t\rightarrow
\varepsilon^{n}$, both amplitude and phase flow equations (\ref{Amp}) and
(\ref{Phase}) now yield linear flow of amplitude and phase relative to the
respective small scale slow, almost constant variables $\eta_{1}$ and
$\eta_{2}$, satisfying $|\eta_{i}^{2}|<<1$ and $|\ddot{\eta}_{i}|<<|\dot{\eta
}_{i}^{2}|<<1$, those vary in a rhythmic manner ( Appendix Sec.B). The exact
(experimental) limit cycle could now be approximated with any desired accuracy
by smooth matching of elementary straight line segments of the form (i)
$z-Z_{0}=k(y-Y_{0})$ and circular arcs of the form (ii) $(y-Y_{0}%
)^{2}+(z-Z_{0})^{2}=a^{2}(\tau_{1})$ over judiciously chosen intervals in $y$
in the $(y,z)$ plane, where $z=\dot{y}$, $k=-\tan(\varepsilon^{n}+\psi
(\tau_{2}))$, $Y_{0}=Y$ and $Z_{0}=\sum_{i}\dot{\tau}_{i}\frac{\partial
\tilde Y}{\partial\tau_{i}}$. Because of the availability of cooperatively evolving
resource of nonlinear scales, such a matching is always possible
theoretically. In the alternative derivation of limit cycle equations in
Appendix Sec. B, we have outlined an approach to gain more analytic
understanding of the rhythmic, cooperative variations of the nonlinear scaling
functions. We hope to address the question of determining the precise analytic
properties of the scaling functions $\phi_{i}(\tilde{\tau})$ corresponding to
a given nonlinear system in future communications.

\begin{figure}[h]
\begin{center}%
\begin{tabular}
[c]{cc}%
\includegraphics[width=3in]{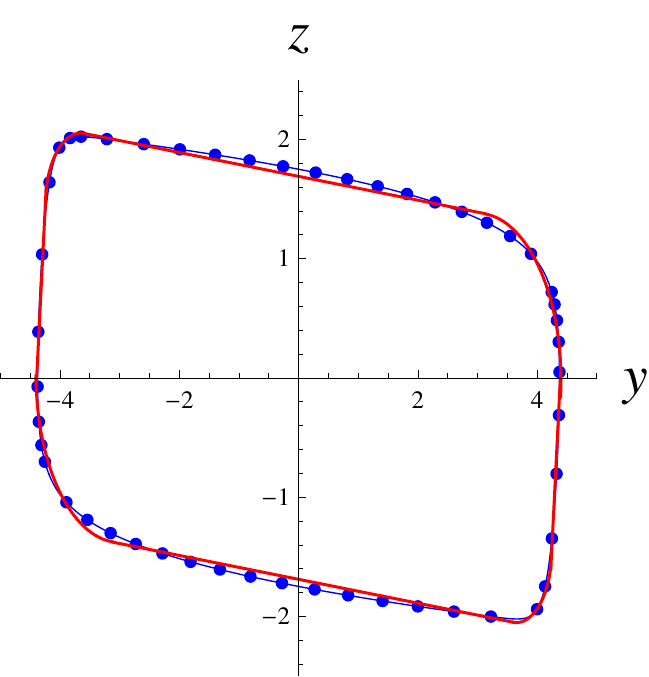} & \includegraphics[
width=3in]{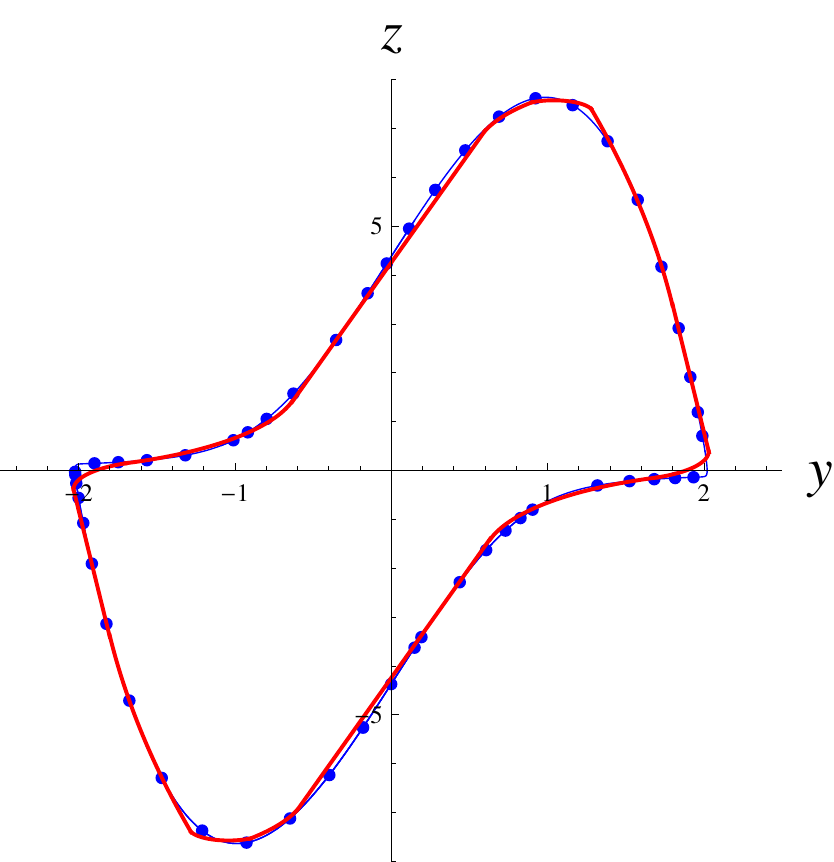}\\
$\left(  a\right)  $ & $\left(  b\right)  $%
\end{tabular}
\end{center}
\caption{Approximate limit cycle for $\left(  a\right)  $ Rayleigh equation
and $\left(  b\right)  $ Van der Pol equation,\newline solid $\left(
\text{red}\right)  $ line for approximate curve, dotted $\left(
\text{blue}\right)  $ line for exact curve $\left(  \varepsilon=5\right)  $.}%
\label{Fig_APP_LC}%
\end{figure}

In Fig. \ref{Fig_APP_LC}$\left(  a\right)  $ and Fig. \ref{Fig_APP_LC}$\left(
b\right)  $, we display the $(y,z)$ phase plane relaxation oscillation for the
Rayleigh and VdP equations with $\varepsilon=5$. In Appendix Sec. C, the
piece-wise smooth matching curves approximating  these cycles are presented in
tabular forms. However, the smoothness at the joining points are achieved at
the level of one decimal only. More accurate approximation may be achieved
with smarter efforts. Judicious choice of slowly varying centres
$(Y_{0},Z_{0})$ and radii $a(\tau_{1})$ of circular arcs of right sizes
(a straight line segment being an arc with sufficiently large radius) should 
give better approximations with a given exact (experimental) 
cycle that can be obtained on a
symbolic computation platform, Mathematica for instance. The phase plane
dynamics of these slowly varying centres and radii is expected to reveal
interesting new insights into asymptotic properties of the nonlinear
oscillation . It transpires from  Appendix Sec. C that  radii, for instance, 
 vary much faster in Rayleigh than that in VdP oscillator, 
in which case radii fluctuate between small and large values through intermediate steps.
This might be compared with fast and slow energy build ups in Rayleigh and VdP relaxation 
oscillations respectively. 
One would like to interpret this phase plane dynamics as
cooperative evolution of multiple nonlinear scales driving amplitude and phase
of the nonlinear oscillation to flow in such a fashion as to generate little
(elementary) circular arcs and linear segments which join smoothly together to form the
complete orbit. Making an accurate plot then boils down to finding right kind
of such arcs and line segments. Intricate dependence of the trajectory itself
into the definition of the nonlinear scales (Appendix Remark 3) tells in
retrospect that one needs to look for a novel iteration scheme that would
allow one to extract the trajectory systematically as a limit process. This
problem will be considered in detail elsewhere.

\section{Concluding Remarks\label{Sec Conclusion}}

In this paper we have presented a comparative study of the homotopy analysis
method and the Renormalization Group method. The approximate formulae for the
amplitudes of the limit cycles of the Rayleigh and the Van der Pol systems are
derived using both the methods and are compared with the exact results. It
turns out that the higher order perturbative calculations based on the
conventional Renormalization group method would fail to give efficient formula
for the limit cycle amplitudes for these nonlinear oscillators. However, an
improved version of the Renormalization group analysis exploiting a novel
concept of nonlinear time is shown to yield efficient amplitude formulae for
all values of $\varepsilon$. Exploiting multiple nonlinear scales of the
associated nonlinear time the improved RG method is also found to yield good
plots for relaxation oscillation orbits for the Rayleigh and VdP systems.

In Appendix we have presented brief review of the nonlinear time formalism and
also given an alternative approach in deriving non-perturbative flow equations
of amplitude and phase of a limit cycle problem. Non-perturbative information
of asymptotic quantities get naturally encoded into nonlinear scales, that can
be exploited judiciously to extract desired asymptotic properties of a
relevant dynamical quantity. More detailed analysis of the nonlinear time
formalism in several other nonlinear systems will be considered in future.

\section*{Acknowledgements}
Authors thank the anonymous  reviewers for constructive criticisms improving the 
quality of the paper.

\section*{Appendix}

\subsection*{A. Formal Structure}

The idea of nonlinear time can be given a rigorous meaning in a nonclassical
extension of the ordinary analysis \cite{DPD1,DPD New}. Recall that the real
number system $\mathbf{R}$ is generally constructed as the metric completion
of the rational field $\mathbf{Q}$ under the Euclidean metric $|x-y|,
\ x,y\in\mathbf{Q}$. More specifically, let S be the set of all Cauchy
sequences $\{x_{n}\}$ of rational numbers $x_{n}\in\mathbf{Q}$. Then $S$ is a
ring under standard component-wise addition and multiplication of two rational
sequences. Then the real number field $\mathbf{R}$ is the quotient space $S/
S_{0}$, where the set $S_{0}$ is the set of all Cauchy sequences converging to
$0\in Q$ and is a maximal ideal in the ring $S$. Alternatively, $\mathbf{R}$
can be considered as the set $[S]$ of equivalence classes, when two sequences
in $S$ are said to be equivalent if their difference belongs to $S_{0}$.

The nonclassical extension $\mathbf{R}^{*}$ of $\mathbf{R}$ is based on a
\emph{finer} equivalence relation that is defined in $S_{0}$ as follows: let
$\{a_{n}\}\in S_{0}$. Consider an associated family of Cauchy sequences of the
form $S_{0a}:=\{A^{\pm}|A^{\pm}=\{a_{n}\times a_{n}^{\pm a^{\pm}_{m_{n}}}\}\}$
where $a^{\pm}_{m_{n}}\neq0$ is Cauchy for $m_{n}>N$ and $N$ sufficiently
large. Clearly, $S_{0a}\subset S_{0}$, and sequences of $S_{0a}$ also
converges to 0 in the metric $|.|$. As $a$ parametrizes sequences in $S$, it
follows that $\underset{\{a\}}\bigcup S_{0a}=S$. Assume further that
$a_{m_{n}}^{\pm}$ respect the \emph{duality structure} defined by $(a_{m_{n}%
}^{-})^{-1}\propto a_{m_{n}}^{+}$ for $m_{n}>N$. The duality structure extends
also over the limit elements: viz., $\mathbf{R}\ni\ (a^{-})^{-1}\propto a^{+}$
where $a^{\pm}_{m_{n}}\rightarrow a^{\pm}$ as $m_{n}\rightarrow\infty$ such
that $a^{\pm}$ are close to 1 in $\mathbf{R}$.

Next define an equivalence relation in $S_{0a}$ declaring two sequences
$A_{1},A_{2}$ in the set $S_{0a}$ equivalent if the associated exponentiated
sequences $a_{m_{n}}^{1}$ and $a_{m_{n}}^{2}$ differ by an element of $S_{0}$
for $m_{n}>N$. In particular, one may impose the condition that $A_{1}\equiv
A_{2}$ if and only if $\exists M \ \mathrm{such \ that} \ a_{m_{n}}%
^{1}=a_{m_{n}}^{2} \ \forall\ m_{n}>M$. Clearly, the usual metric $|.|$ fails
to distinguish elements belonging to two distinct such finer equivalent
classes. However, the metric defined as the natural logarithmic extension of
the Euclidean norm, generically called the \emph{asymptotically visibility
metric} is introduced by $h(A_{1},A_{2}) =\underset{n\rightarrow\infty}%
\lim|\log_{|A_{0}|^{-1}} |A_{1}-A_{2}|/|A_0||$ where $A_{0}=\{a_{n}\}\in S_{0}$. The
sequence $A_{0}$ is said to define a natural scale relative to which elements
in $S_{0}$ gets nontrivial values and hence become distinguishable. The limit
exists because of concerned sequences $a_{m_{n}}^{\pm}$ being Cauchy. Note
that the mapping $h:S_{0}\rightarrow\mathbf{R^{+}}$ defined by $h(A)=\underset
{n\rightarrow\infty}\lim|\log_{|A_{0}|^{-1}} |A|/|A_{0}||$ is actually a
nontrivial norm \cite{DPD1,DPD New} (for simplicity of notation, we use same
symbol to denote both the norm and metric).

The extended real number system $\mathbf{R}^{*}$ admitting duality induced
\emph{fine structure} is given, by definition, as the equivalence class under
this finer equivalence relation viz., $\mathbf{R}^{*}:=S/S_{0}$ when
convergence is induced naturally by the asymptotically visibility metric
$h(x,y)$. Clearly, under the usual norm $|.|$, $\mathbf{R}^{*}$ reduces to
$\mathbf{R}$ as the exponentiated elements $a^{\pm}$ are essentially
invisible. The natural application of the visibility norm on $\mathbf{R}^{*}$
is activated in the following steps. For any two distinct elements $x,\ y\in
{\bf R}\subset\mathbf{R}^{*}$, \emph{set}, by definition, $h(x,y)=0,\ x\neq y$;
$h(x,y)$ being nontrivial only for $y\in x+S_{0}$. This choice is natural as
for any element $x\in\mathbf{R}$, the corresponding limiting $h$ norm viz,
$h(x)=h(0,x)=\underset{n\rightarrow\infty}\lim\log_{|A_{0}|^{-1}} |x/A_{0}|$=1 and
so, the definition $h(x,y)=0, \ \forall x,y\in\mathbf{R}$ makes sense. 
For nontrivial values of $h(x,y),
\ x,\ y\in\mathbf{R}^{*}$, the definition of the visibility metric extends
over to $h(x,y)=\underset{\varepsilon\rightarrow0}\lim|\log_{\varepsilon^{-1}}
|\frac{x-y}{\varepsilon}||$, which exists by construction.

Next, consider the metric $d: \mathbf{R}^{*}\rightarrow {\bf R}^{+}$ by
$d(x,y)=|x-y|+h(x,y)$. Clearly, ${d}(x,y)=|x-y|$ for any $x,\ y\in {\bf R}$ and
${d}(x,y)=h(x,y)$ for $x,\ y\in\mathbf{R}^{*}-{\bf R} $ and hence $(\mathbf{R}%
^{*},{d})$ is a complete metric space. The metric $h(x,y)$ acting nontrivially
on $S_{0}$ is essentially an ultrametric: $h(x,y)\leq\max\{h(x,z),h(z,y)\}$.
This follows immediately from the observation that $h$ maps $\mathbf{R}$ to
the singleton set $\{1\}$. Further, the ultrametric $h$ must be discretely
valued \cite{DPD1} and hence the nontrivial value set of $h$ viz., $h(S_{0})$
is countable. As a consequence, the set $S_{0}$ is totally disconnected and
perfect in the induced topology.

More detailed analytic aspects (including the idea of smooth jump
differentiability and jump derivative) of the extended system $\mathbf{R}^{*}$
equipped with the metric $d$ will be reported elsewhere \cite{DPD New}. Here,
we make a few relevant remarks.

1. Even as the size of a $\delta-$ neighbourhood of a point $x\in\mathbf{R}$
vanishes linearly, the same for $x^{*}\in\mathbf{R^{*}}$ need not vanish at
the same rate and may only vanish at a slower rate $\delta h(\delta)$. The
real number model $\mathbf{R}$ is called the hard or string model when the
space $\mathbf{R^{*}}$ is called the soft or fluid model of real numbers
\cite{DPD Sen}. The ordinary differential measure $dx$ gets extended in
$\mathbf{R^{*}}$ as $d(h(x)x)$.

2. Consider the open interval $(\delta,\delta^{-1})\subset\mathbf{R^{*}}$. In
the asymptotic limit $\delta\rightarrow0^{+}$, the duality structure
identifies the right neighbourhood of $\delta$ with the left neighbourhood of
$\delta^{-1}$ in a nontrivial manner. As a consequence, the linear
(translation) group action on $\mathbf{R}$ is extended to a nonlinear
SL($2,R$) group on $\mathbf{R^{*}}$. Infact, the translation subgroup acts on
$\mathbf{R}$, when the inversion acts nontrivially only on $\mathbf{R^{*}}$ in
the sense that the visibility norm $h$ is invariant under inversion $\hat i:
\ h(\hat i A)=\hat i(h(A))$ where $\hat i(A)=\{a_{n}^{-1}\times a_{n}%
^{(a^{-}_{m_{n}})^{-1}}\}, \ A=\{a_{n}\times a_{n}^{-a^{-}_{m_{n}}}\}$. For a
translation $T_{r}$ by a shift $r$, on the other hand, $h(T_{r}(A))=h(A)$ and
hence $T_{r}(A)=A \ \Rightarrow r=0$ (i.e. $T$ acts trivially on $\bf{R^*-R}$). Above two
\emph{salient} properties of the duality structure are expected to have
significant application in nonlinear problems.

3. To give an example of the intricate nonlinear structure that can get
encoded into a well behaved (smooth) function in $\mathbf{R}$, let us consider
the simplest case of a real variable $x$. In $\mathbf{R^{*}}$ the variable $x$
gets extended to, say, $X=x e^{\phi(\log X)}$. The function $\phi$ exposing
the nonlinear dependence is also assumed to be differentiable. Differentiating
$X$ with respect to $x$ one gets $xX^{\prime}(1-\phi^{\prime})=X$, where
$\prime$ denotes derivation with the argument. We now assume that $\phi(\log
X)$ is vanishingly small (i.e. less than accuracy level $\delta$ in any given
application) for $0<x<\infty$ and O(1) when $\log X>>1$ i.e. $x\rightarrow0$
or $\infty$. As a consequence, existence of $\phi$ is felt only in the
asymptotic neighbourhoods (Remark 2) of 0 or $\infty$. We now make a further
assumption that $\phi^{\prime}=0$ almost everywhere in an asymptotic
neighbourhood, but everywhere in $0<x<\infty$. Then $X$ satisfies
$xX^{\prime}=X$ a.e. in $\mathbf{R^{*}}$. Thus ordinary variable
$x\in\mathbf{R}$ gets extended in $\mathbf{R^{*}}$ as $X$ which has the
intermittent property of a Cantor devil's Stairecase function in an asymptotic
neighbourhood. Since, under duality structure, such a neighbourhood has
ultrametric topology, $X$ in fact satisfies the above scale invariant equation
everywhere in $\mathbf{R^{*}}$, because ordinary non-differentiability at the
points of the associated Cantor set is removed by inversion mediated jump
increments \cite{DPD1, DPD New}. This example tells that an ordinary function
can have nonlinear and nonlocal functional dependence with itself, along with
rhythmic (intermittent) variability that can have significant amplification in
an asymptotic sector.

4. The asymptotic scaling variables $h_{0}(\varepsilon)$ and $H_{RG}%
(\tilde\tau)$ introduced in Sec. 4 correspond to the associated visibility
norm $h(A)$ defined above. A real variable $t\in\mathbf{R}$ approaching
asymptotically either to 0 or $\infty$ has natural images in $\mathbf{R}^{*}$
in the form $\tau_{0}=t\times t^{-h^{-}(\varepsilon t)}, \ h^{-}(\varepsilon
t)<1$ and $\tau_{\infty}=t\times t^{h^{+}(\varepsilon t)}, \ h^{+}(\varepsilon
t)>1$ respectively. The scaling exponents $h^{\pm}$ encode asymptotic scaling
information of a given nonlinear system. Further, $(h^{-}(\varepsilon
t))^{-1}\propto h^{+}(\varepsilon t)$ by duality. In Sec.4, we discuss how
such information can be systematically extracted in the case of a limit cycle
for a nonlinear oscillator (see also below).

5. The fine structures in $\mathbf{R^{*}}$ remain inactive (passive/hidden) in
absence of any stimulus, either intrinsic or external. In presence of an
external input, say, the actions of the nontrivial component of the metric $d$
and the associated duality structure become manifest. The RG analysis makes
a room for direct implementation of the intrinsically realized duality structure
in the context of a nonlinear system in the soft model $\mathbf{R^{*}}$.

\subsection*{B. Application: Limit Cycle}

Consider a general nonlinear oscillator given by
\begin{equation}
\label{NO}\ddot{x}+x= \varepsilon f(x,\dot x)
\end{equation}
We assume $f$ such that the system admits a unique isolated cycle for
$\varepsilon>0$ and other relevant parameter values. For a finite nonlinearity
$\varepsilon>1$ (say), the usual linear time $t$ is extended to one enjoying
\emph{right} asymptotic correction $t\rightarrow\ T_{i}=t\phi_{i}(\tilde
\tau(t))$ as $t\rightarrow\infty$ through linear scales $\varepsilon^{i}$,
where $\phi_{i}(\tilde\tau)$ stands succinctly for the nontrivial
intrinsically generated \emph{slowly varying} scaling components arising from
the associated visibility norm. Here, $\tilde\tau$, as usual denotes an O(1)
rescaled variable in the neighbourhood of linear scales $\varepsilon^{i}$. In
the case of a nonlinear planar autonomous system the relevant dynamical
quantities are only amplitude and phase of the nonlinear oscillation and so we
have only two asymptotic scaling functions $\phi_{i}(\tilde\tau), \ i=1,2$
which get \emph{selected} naturally so as to facilitate direct
non-perturbative calculation of the asymptotic properties i.e. the amplitude
and phase of the limit cycle of the system. An implementation of this
non-pertubative scheme in the perturbative RG formalism is presented in
Sec.4.1 for computation of the amplitude of the concerned oscillators. In
Sec.4.2, we have demonstrated that the computed plot of the limit cycle for
the Rayleigh and VdP oscillators could be matched arbitrarily closely for
appropriate choices of the slowly varying nonlinear time when the amplitude
and phase of the unperturbed periodic solution \emph{flow linearly} in the
appropriately chosen nonlinear scaling time variables.

Here, we give an alternative derivation of the nonperturbative
\emph{relaxation oscillation} flow equations ab-initio from the slowly varying
nonlinear time in the context of the Rayleigh equation (\ref{Rayleigh}) with
$\varepsilon>>1$. It will transpire that the new approach is free of any
divergence problem because of its inbuilt RG cancellations via duality
principle. Since we are interested in the planar limit cycle properties, we
assume that all the relevant quantities e.g. the solution $y$, amplitude $a$
and phase $\psi$ are functions of asymptotic time variable $t\sim
\varepsilon^{n},\ n>>1$ and the associated nontrivial scaling variables
$\tau_{1}=\phi_{1}(\tilde{\tau})$ and $\tau_{2}=\phi_{2}(\tilde{\tau})$ for a
rescaled $\tilde{\tau}\sim O(1)$. Higher order scaling variables $\tau
_{n},\ n>2$ of the nonlinear structure of time variable may become relevant
for a non-planar system. Accordingly, we write the ansatz $y(t,\tau_{1}%
,\tau_{2})=y_{0}(t)+Y_{1}(\tau_{1})+Y_{2}(\tau_{2})$ for the limit cycle
solution involving multiple time scales (only three for the planar system).
Assuming slow variations of nonlinear scales $\tau_{i}=\phi_{i},\ i=1,2$, viz.
$|\phi_{i}^{\prime\prime}|<<|{\phi_{i}^{\prime}}^{2}|<<1,\ \phi_{i}^{\prime
}=\frac{d\phi_{i}}{d\tilde{\tau}}$, as $t\rightarrow\varepsilon^{n}%
,\ n\rightarrow\infty$ and noting that $\frac{dy}{dt}=\frac{\partial y_{0}%
}{\partial t}+\sum_{i}\dot{\tilde{\tau}}\phi_{i}^{\prime}{\frac{\partial
Y_{i}}{\partial\tau_{i}}}$ etc., the Rayleigh equation (\ref{Rayleigh})
simplifies to
\begin{equation}
\frac{\partial^{2}y_{0}}{\partial t^{2}}+y_{0}+\sum_{i}Y_{i}=\varepsilon
(\frac{\partial y_{0}}{\partial t}+\sum_{i}\dot{\tilde{\tau}}\phi_{i}^{\prime
}{\frac{\partial Y_{i}}{\partial\tau_{i}}})-\frac{\varepsilon}{3}%
((\frac{\partial y_{0}}{\partial t})^{3}+3(\frac{\partial y_{0}}{\partial
t})^{2}\sum_{i}\dot{\tilde{\tau}}\phi_{i}^{\prime}{\frac{\partial Y_{i}%
}{\partial\tau_{i}}})\label{mulscale}%
\end{equation}
where we drop all higher order terms involving $\phi_{i}^{\prime\prime}$ and
${\phi_{i}^{\prime}}^{2}$. Assuming $y_{0}(t)=a(\tau_{1},\tau_{2})\cos
(t+\psi(\tau_{1},\tau_{2}))$ with \emph{flowing} amplitude and phase in
scaling times $\tau_{1}$ and $\tau_{2}$ so that
\begin{equation}
\frac{\partial^{2}y_{0}}{\partial t^{2}}+y_{0}=0\label{linear1}%
\end{equation}
we next get a simplified linearized evolution for the nonlinear components of
the asymptotic limit cycle solution in the form
\begin{equation}
(1-(\frac{\partial y_{0}}{\partial t})^{2})\sum\dot{\phi}_{i}\frac{\partial
Y_{i}}{\partial\tau_{i}}=\{\frac{1}{3}(\frac{\partial y_{0}}{\partial t}%
)^{3}-\frac{\partial y_{0}}{\partial t}\}+\varepsilon^{-1}\sum Y_{i}%
\label{linear2}%
\end{equation}
where $\dot{\phi}_{i}=\frac{d\phi_{i}}{dt}$. As a consequence, under the
assumption of slowly varying nonlinear time scales, a second order nonlinear
planar system (\ref{Rayleigh}) would decompose into a linear second order
partial differential equation (\ref{linear1}) for the zero level solution
$y_{0}$ and an associated first order partial differential equation
(\ref{linear2}) for the nonlinear scale dependent components $Y_{i}$. Clearly,
analogous decomposition holds actually for a larger class of planar autonomous
systems (\ref{NO}) having a unique limit cycle solution. Extension of this
result to multiple limit cycles would be considered separately.

We note here that since the system (\ref{linear1}) and (\ref{linear2}) is
\emph{under determined}, there is room for further restrictions to solve the
system self-consistently. To re-derive the RG flow equations (\ref{Amp}) and
(\ref{Phase}) from (\ref{linear2}), we now make following assumptions: we
write (i) $a(\tau_{1},\tau_{2})=a(\tau_{1}),\ \psi(\tau_{1},\tau_{2}%
)=\psi(\tau_{2})$ so that amplitude varies slowly with first order scale
$\tau_{1}$ when the second order scale $\tau_{2}$ and phase $\psi$ remain
almost constant. On the other hand as $a$ stabilizes to an almost constant
value, the phase begins to flow, though slowly with the second order scale
$\tau_{2}$. Such slow, almost constant, rhythmic cooperative variations of
$\tau_{1}$ and $\tau_{2}$ are modeled, depending on the specific problem under
consideration (see below and c.f. Sec. 4.2), to retrieve the RG flow equations
correctly. As shown in the example (Remark 3, Appendix Sec.A) such a rhythmic
nonlinear variation does exist in an ultrametric neighbourhood in
$\mathbf{R^{\ast}}$. To further quantify the slow variation of dynamical
variables, we next impose the condition that (ii) \emph{the total variation of
the exact solution $y(t,\tau_{1},\tau_{2})$ with respect to each slow variable
$\tau_{i}$ along the full periodic cycle $C$ must vanish viz, $\int_{C}%
\frac{\partial y}{\partial\tau_{i}}dt=0$ for each $i$}. To avoid trivialities
i.e. $\int_{C}\cos(t+\psi)dt=0$ etc., we, however, evaluate the concerned
integrals only on the quarter cycle, with the understanding that phase shifts
of $\pi/2$ are absorbed in the definition of $\psi$.

In the sufficiently large $\varepsilon>1$ relaxation oscillation, one can
further simplify (\ref{linear2}) by dropping the $\varepsilon^{-1}$ term to
obtain
\begin{equation}
\label{linear3}\sum\dot\phi_{i}\frac{\partial Y_{i}}{\partial\tau_{i}}%
=\frac{\frac{1}{3}(\frac{\partial y_{0}}{\partial t})^{3}-\frac{\partial
y_{0}}{\partial t}}{1-(\frac{\partial y_{0}}{\partial t})^{2}}\equiv
\Phi(y_{0t}), \ y_{0t}=\frac{\partial y_{0}}{\partial t}%
\end{equation}
To make contact with RG flow equations (\ref{Amp}) and (\ref{Phase}) one now
exploits \emph{the freedom of right choice} in the functional forms of
nonlinear scales. For the Rayleigh equation, we now set for slow,
cooperatively active functional dependence (a) $\dot\phi_{1}=\Phi(y_{0t}%
)S_{1}^{-1}(a,\psi,t), \ \dot\phi_{2}=0$ and (b) $\dot\phi_{1}=0, \ \dot
\phi_{2}=\Phi(y_{0t})S_{2}^{-1}(a,\psi,t)$ for successive slow variations, as
described in (i), of the scales $\tau_{1}$ and $\tau_{2}$ respectively, where,
$S_{1}= \frac{1}{2}a(\frac{a^{2}}{4}-1)\cos(t+\psi)$ and $S_{2}=\frac{1}%
{8}(1-\frac{a^{4}}{32})\sin(t+\psi)$ (recall the example in Remark 3 of Sec.A
above highlighting wide possible choices and intricate functional dependence).
These choices for $S_{1}$ and $S_{2}$ would yield the RG flow equations when
condition (ii) is invoked.

Note that the relations in both (a) and (b) are truly nonlinear; the dynamical
variables $a$ and $\psi$ in $S_{1}$ and $S_{2}$ depend implicitly in $\phi
_{1}$ and $\phi_{2}$ respectively, which, in turn are \emph{slowly varying} as
the linear parameter $t$ is assumed to vary in a neighbourhood of
$\varepsilon^{n}$ for a large but fixed $n$. Invoking the global slow
variation condition (ii) for each $i$, in conjunction with the ansatz (a) and
(b), one finally deduce the amplitude and phase flow equations (\ref{Amp}) and
(\ref{Phase}) in slow variables $\tau_{1}$ and $\tau_{2}$ respectively.

In the present format, the flow equations, however, have got \emph{new}
interpretations: Amplitude and phase must flow in successive rhythmic manner;
phase remains almost constant (i.e. $\frac{\partial\psi}{\partial\tau_{i}}=0$
for each $i$) when amplitude varies slowly with $\tau_{1}$ towards an almost
constant value. Subsequently, the flowing of $a$ is halted temporarily (i.e.
$\frac{\partial a}{\partial\tau_{i}}=0$), initiating flowing of $\psi$ in next
level variable $\tau_{2}$. This rhythmic oscillation would obviously continue
indefinitely over a cycle. The RG flow equations could be treated as
non-perturbative because of \emph{implicit} connections of nonlinear scaling
time functions with amplitude and phase via intrinsically defined duality
principle (c.f. Remark 3 above).

To summarize, based on perturbative RG formalism we have presented an
alternative approach in deriving non-perturbative flow equations of relevant
asymptotic dynamical quantities of a planar autonomous limit cycle problem,
from nonlinear scale invariant time scales, which become available in an
extended analytic framework incorporating duality structure. Non-perturbative
information of asymptotic quantities get naturally encoded into nonlinear
scales, that can be exploited judiciously to extract desired asymptotic
analytic properties of a relevant dynamical quantity. An algorithmic procedure
of extracting such information is explained in estimating both the limit cycle
amplitude and trajectory for Rayleigh and VdP equations.

\subsection*{C. Matching Arcs and Segments}

$\left(  a\right)  $ Piecewise smooth matching curves for upper half of the
approximate Rayleigh limit cycle for $\varepsilon=5$:

$z\left(  y\right)  =\left\{
\begin{array}
[c]{cc}%
0.02-\sqrt{1.96-\left(  y+3\right)  ^{2}} & -4.96<y\leq-4.393\\
10y+43.9 & -4.393<y\leq-4.23\\
1.46+\sqrt{0.35-\left(  y+3.65\right)  ^{2}} & -4.23<y\leq-3.6\\
-0.1y+1.689 & -3.6<y\leq3.12\\
-0.02+\sqrt{1.96-\left(  y-3\right)  ^{2}} & 3.12<y\leq4.96
\end{array}
\right.  $

$\left(  b\right)  $ Piecewise smooth matching curves for upper half of the
approximate VdP limit cycle for $\varepsilon=5$:

$z\left(  y\right)  =\left\{
\begin{array}
[c]{cc}%
-0.388+\sqrt{0.352-\left(  y+1.438\right)  ^{2}} & -2.05<y\leq-1.7\\
1.8-\sqrt{3.2-\left(  y+2.38\right)  ^{2}} & -1.7<y\leq-0.633\\
4.5y+4.25 & -0.633<y\leq0.6\\
6.5+\sqrt{2.76-\left(  y-2.2\right)  ^{2}} & 0.6<y\leq0.9\\
7.325+\sqrt{0.063-\left(  y-1.04\right)  ^{2}} & 0.90<y\leq1.28\\
0.38+\sqrt{1530-\left(  y+37.2\right)  ^{2}} & 1.28<y\leq1.8\\
-13y+26.8 & 1.8<y\leq2.033\\
0.388+\sqrt{0.352-\left(  y-1.438\right)  ^{2}} & 2.033<y\leq2.05
\end{array}
\right.  $

\end{document}